\documentclass[]{article}

\usepackage{graphicx} 
\usepackage{float}
\usepackage{booktabs}

\usepackage{epstopdf}
\usepackage{amsfonts}
\usepackage{amsmath}
\usepackage{amssymb}
\usepackage{fancyhdr}
\usepackage{titlesec}
\usepackage{indentfirst}
\usepackage{booktabs}
\usepackage{verbatim}
\usepackage{color}
\usepackage{xcolor}
\usepackage{caption}
\usepackage{amsthm}
\usepackage{mathtools}
\usepackage{siunitx}
\usepackage{placeins}  

\newcommand{\steptitle}[1]{\FloatBarrier\noindent\textbf{#1}\par\vspace{0.4em}}

\usepackage[normalem]{ulem}
\usepackage{subcaption} 

\usepackage{amsthm,amsmath,amssymb,dsfont}
\usepackage[colorlinks]{hyperref}

\usepackage[page,header]{appendix}
\usepackage{titletoc}

\DeclareCaptionLabelFormat{redfigure}{\textcolor{red}{#1~#2}}

\numberwithin{equation}{section}
\newtheorem{theorem}{Theorem}[section]

\newtheorem{remark}[theorem]{Remark}

\newcommand{\rd}{\mathrm{d}}

\newcommand{\im}{\mathrm{i}}
\newcommand{\e}{\mathrm{e}}
\newcommand{\mQ}{\mathcal{Q}}
\newcommand{\mP}{\mathcal{P}}
\newcommand{\M}{\mathcal{M}}

\newcommand{\G}{\mathcal{G}}

\newcommand{\Sb}{\mathbb{S}}
\newcommand{\Sd}{\mathbb{S}}
\newcommand{\R}{\mathbb{R}}

\begin{document}

\title{Solving the BGK model and Boltzmann equation by Fourier Neural Operator with conservative constraints}

\author{Boyun Hu and Kunlun Qi}

\author{
	Boyun Hu\footnote{School of Mathematics, University of Minnesota - Twin Cities, Minneapolis, MN 55455, USA (hu000608@umn.edu ).} \  \ and \
	Kunlun Qi\footnote{Department of Computational Mathematics, Science and Engineering and Department of Mathematics, Michigan State University, East Lansing, MI 48824, USA (qikunlun@msu.edu, kunlunqi.math@gmail.com).}
	} 

\date{}

\maketitle

\begin{abstract}
The numerical approximation of the Boltzmann collision operator presents significant challenges because of its high dimensionality, quadratic nonlinearity, and nonlocal integral structure. In this work, we develop a Fourier Neural Operator (FNO)-based framework for learning finite time evolution operators associated with the spatially homogeneous Boltzmann equation and its simplified BGK model in different velocity space dimensions. The framework learns mappings between distribution functions in either sequence-to-sequence or point-to-point form and, once trained, enables rapid evaluation without repeated classical time evolution. To improve physical consistency, we introduce a conservative FNO (C-FNO) by augmenting the training loss with soft penalties for deviations in mass, momentum, and energy from the conservation law.
Several numerical experiments are presented to demonstrate that the learned operators accurately reproduce the principal collision dynamics, while C-FNO generally reduces macroscopic drift relative to the standard FNO, highlighting its potential as a promising framework for physics-constrained operator learning in kinetic theory and other nonlinear integro-differential equations.
\end{abstract}

{\small 
Fast Fourier Transform, Conservation Law.

}

\tableofcontents

\section{Introduction}

The Boltzmann equation, introduced by Maxwell and Boltzmann, is a fundamental model in kinetic theory that describes the time evolution of the particle distribution function in rarefied gases \cite{CC, Cercignani}. Unlike macroscopic fluid models, it remains valid in regimes where the continuum assumption fails \cite{Villani02}. In particular, when the mean free path of particles is comparable to the characteristic length scale of the system, a kinetic description becomes essential. Such conditions arise in rarefied gas dynamics as well as in microscale and nanoscale devices, where characteristic dimensions are small.
Beyond its broad range of applications in science and engineering, the Boltzmann equation plays a central role in mathematics and physics, underpinning the statistical foundations of non-equilibrium thermodynamics.

Albeit its fundamental importance in the statistical physics and other applications, solving the Boltzmann equation numerically remains a longstanding challenge \cite{pareschi}. The main trouble-maker is the complicated collision operator, posed in a high-dimensional integral form, which leads to severe computational complexity often referred to as the curse of dimensionality. In addition, the highly nonlinear and nonlocal structure of the collision operator demands expensive computational cost to achieve accuracy.
To mitigate the curse of dimensionality, particle-based methods such as the direct simulation Monte Carlo (DSMC) scheme \cite{Nanbu80,Bird} have been widely employed, offering notable advantages in efficiency and ease of implementation. However, these stochastic type methods often suffer from slow convergence and become particularly inefficient in non-steady or low-speed regimes. In contrast, deterministic approaches have seen rapid progress in recent decades, driven by advances in computing power \cite{pareschi}. Among them, the Fourier–Galerkin spectral method provides a powerful framework for approximating the Boltzmann collision operator \cite{PP96, BR99}. This approach not only achieves spectral accuracy but also admits efficient extensions through accelerated algorithms \cite{MP06, MPR13, WWSRZ2013, GHHH17, HQ2020, HQY2021, QSW2025}. We also refer to \cite{CFY2018, WC2019, FPR15, AGT16, PR2022} for other types of spectral method in solving the Boltzmann equation.

While traditional numerical methods have been successful, they are often computationally expensive, particularly for high-dimensional problems in statistical physics, fluid dynamics, and many-body dynamical systems. To achieve high accuracy, these solvers must discretize space and time into extremely fine grids, requiring the solution of a massive number of equations. This high computational cost represents a significant bottleneck for complex scientific and engineering applications.
Therefore, a recent paradigm shift in scientific computing involves the development that leverage machine learning and neural network. For architectures that do not use data, it is worth mentioning the recent progress in applying neural network in solving the multiscale kinetic equation, for instance, inspired by the PINN framework \cite{RPK2019}, Jin-Ma-Wu in \cite{JMW2023, JMW2024} propose novel Asymptotic-Preserving Neural Networks (APNNs) for tackling multiscale time-dependent kinetic problems, encompassing the linear transport equation and BGK model and covering all ranges of Knudsen number. 
Furthermore, for those that partially use data, Li-Wang-Liu-Wang-Dong in \cite{LWLWD2024} take advantage of the neural networks to reduce complexity via sparse or low-rank representations of the distribution function in velocity space.
For those data-driven methodologies in Boltzmann equation, Xiao-Frank in \cite{XF2021, XF2023} and Miller-Roberts-Bond-Cyr in \cite{MRBC2022} develop different frameworks to blend classical kinetic modeling with neural network surrogates, while enforcing physical structure.
In addition to serving as surrogate models or solving equations, the neural networks are also applied in the moment closure modeling of the kinetic equation, for which we refer to \cite{LDW2023, Huang1, SLFH2025, HLQW2025} and the references therein to the readers for recent progress.

More recently, operator learning has emerged as a novel approach that leverages the underlying structure of complex problems \cite{kovachki2024operator}. Rather than directly solving equations through fine-grained discretization, these methods learn mappings from inputs to solutions using training data. Once trained, such solvers can rapidly generate predictions for new problems, sometimes outperforming traditional numerical schemes in terms of speed, making them especially attractive for computationally intensive applications. 
Building on this idea, You-Xu-Cai \cite{YXC2024} proposed the multi-scale Fourier neural operator (MscaleFNO) to enhance learning in oscillatory function spaces.
Among various operator learning frameworks, the Fourier Neural Operator (FNO) has demonstrated remarkable success \cite{Li-ICLR2021, LHLA2023}. By learning mappings between function spaces in the Fourier domain, the FNO generalizes across spatial resolutions, allowing models trained on coarse data to yield accurate high-resolution predictions. This capability significantly reduces the computational cost of high-fidelity simulations while maintaining predictive accuracy. 
Despite the growing success of operator learning for various classes of partial differential equations, its application to kinetic equations—particularly the Boltzmann equation—remains relatively unexplored. To the best of our knowledge, Liu-Cai-Xu \cite{LCX2023} applied the MscaleDNN framework to the Poisson–Boltzmann equation in complex geometries, while Lee-Jung-Lim-Hwang \cite{LJLH2025} introduced the Fourier Neural Spectral Network (FourierSpecNet), a hybrid approach combining the Fourier spectral method with deep learning to efficiently approximate the Boltzmann collision operator in Fourier space.

Considering the concerns of the traditional numerical methods in solving the Boltzmann equation and the development of the operator learning framework, in this work, we aim to advance the application of operator learning in kinetic theory by developing a Fourier Neural Operator (FNO) framework for learning the evolution operator of homogeneous Boltzmann equation and its simplified counterpart, the Bhatnagar–Gross–Krook (BGK) model, across different dimensions. In addition, to enhance the physical fidelity of the learned operator, we further modify the original FNO architecture to incorporate intrinsic physical information, adding the conservation of key macroscopic quantities such as mass, momentum, and energy in the loss function as soft constraints. This physically informed modification not only improves the interpretability and stability of the learned models but also possibly promotes better generalization across a wide range of kinetic regimes. The proposed framework thus provides a promising step toward physics-constrained operator learning for other kinetic equations.

The rest of this paper is organized as follows. In the next section \ref{sec:boltzmann}, we give a basic introduction of Boltzmann equation and its simplified BGK model.
Our proposed methodology and data generation process are then described in section \ref{sec:method}. 
In section \ref{sec:numerical}, extensive numerical experiments are presented to demonstrate the eﬃciency and accuracy of the proposed framework. 
The paper is finally concluded in section \ref{sec:conclusion}.

\section{The BGK model and Boltzmann equation}
\label{sec:boltzmann}
The kinetic equation reads
\begin{equation}
\partial_{t} f + v \cdot \nabla_x f=\mQ[f], \quad t>0, \ x\in \Omega\subset\mathbb{R}^d, \ v\in \mathbb{R}^d,
\end{equation}
where $f = f(t,x,v)$ is the probability density function of finding particles with velocity $v$, in the small neighborhood of the position $x$ at time $t$, and $\mQ[f]$ is generally defined as an operator to describe the collision between particles, acting on $f$ only through the velocity $v$. Therefore, in this paper, we mainly consider the spatially homogeneous case, where the distribution function $f$ is independent of position $x$,
\begin{equation}\label{kin}
\partial_{t} f = \mQ[f], \quad t>0, \ v \in \mathbb{R}^d,
\end{equation}
and the collision operator $\mQ[f]$ includes the classical Boltzmann collision operator and its simplified BGK collision operator as follows.

\textbf{Boltzmann collision operator}: $\mQ[f]$ is the nonlinear operator describing the binary collisions among particles:
\begin{equation}\label{Bol}
\mQ[f] = \mQ_{\text{Bol}}(f,f)(v) := \frac{1}{\text{Kn}}\int_{\mathbb{R}^d}\int_{\Sb^{d-1}}\mathcal{B}(v-v_*,\sigma)\left[f(v')f(v_*')-f(v)f(v_*)\right]\,\rd{\sigma}\,\rd{v_*}.
\end{equation}
In the formula above, $\text{Kn}$ is the Knudsen number, which takes $\text{Kn}=1$ for the Boltzmann equation throughout this paper, $(v',v_*')$ and $(v,v_*)$ represent the velocity pairs before and after the collision process, which satisfy the conservation of momentum and energy:
\begin{equation*} 
v' + v_{*}' = v + v_{*}, \quad  |v'|^{2} + |v_{*}'|^{2} = |v|^{2} + |v_{*}|^{2},
\end{equation*}
so that $(v',v_*')$ can be expressed in terms of $(v,v_*)$ as
\begin{equation*}  
v'=\frac{v+v_*}{2}+\frac{|v-v_*|}{2}\sigma, \quad  v_*'=\frac{v+v_*}{2}-\frac{|v-v_*|}{2}\sigma,
\end{equation*}
where $\sigma$ is a unit vector varying over the unit sphere $\Sb^{d-1}$. \\
Additionally, in the Boltzmann operator \eqref{Bol}, the collision kernel $\mathcal{B}$ is a non-negative function depending only on the relative velocity $|v-v_*|$ and the cosine of the deviation angle $\theta$ (angle between $v-v_*$ and $v'-v_*'$) \cite{Villani02}. Hence, $\mathcal{B}$ is often written as
\begin{equation} 
\mathcal{B}(v-v_*, \sigma)=B(|v-v_*|,\cos \theta), \quad \cos\theta=\frac{\sigma\cdot (v-v_*)}{|v-v_*|}.
\end{equation}
More specifically, in the case of inverse power law potentials $U(r) = r^{-(\mathrm{s}-1)}, 2< \mathrm{s} < \infty$, where $r$ is the distance between two interacting particles, $B$ can be separated as the kinetic part and angular part:
\begin{equation}\label{Bb}
B(|v-v_*|,\sigma) = b(\cos\theta) \Phi(|v-v_*|), \quad \cos\theta=\frac{\sigma\cdot (v-v_*)}{|v-v_*|},
\end{equation} 
where kinetic collision part $\Phi(|v-v_*|)=|v-v_*|^{\gamma}$, $\gamma = \frac{\mathrm{s}-5}{\mathrm{s}-1}$, includes hard potential $ (\gamma>0) $, Maxwell molecules $ (\gamma =0) $ and soft potential $ (\gamma<0) $. In addition, the angular collision part $b(\cos\theta)$ is an implicitly defined function, asymptotically behaving as, when $\theta \rightarrow 0^{+}$,
\begin{equation}\label{noncutoffnu}
\sin^{d-2} \theta b(\cos\theta) \big|_{\theta\rightarrow 0^{+}} \sim K\theta^{-1-\nu}, \quad \nu = \frac{2}{\mathrm{s}-1}, \quad 0<\nu<2 \quad \text{and} \quad  K >0.
\end{equation}

\textbf{BGK collision operator}: as a simplified model, the Bhatnagar–Gross–Krook (BGK) model \cite{bhatnagar1954model} replaces the complicated integral operator \eqref{Bol} by a simple output term, 
\begin{equation}\label{BGK}
    \mQ[f] = \mQ_{\text{BGK}}(f)(v) := \frac{f_{\text{eq}}(v)-f(v) }{\text{Kn}},
\end{equation}
where $f_{\text{eq}}$ is the equilibrium state, or so-called the Maxwellian distribution,
\begin{equation}\label{feq}                                                                                 
f_{\text{eq}}(v) := \M(v) = \frac{\rho}{\left(2\pi T\right)^{\frac{d}{2}}} \exp\Big( -\frac{|v-u|^{2}}{2T} \Big),
\end{equation}
which depends on the density $ \rho $, bulk velocity $ u$, and temperature $ T $. These macroscopic quantities are derived from the distribution function $f$ in the sense that
\begin{equation}\label{rho-u-t}
\rho = \int_{\mathbb{R}^d} f \,\rd v, \quad m := \rho u=\int_{\mathbb{R}^d} f v \,\rd v, \quad E:= \frac{1}{2}\rho |u|^2 + \frac{d}{2}\rho T = \frac{1}{2} \int_{\mathbb{R}^d} f |v|^{2} \,\rd v,
\end{equation}
where $m$ is the momentum and $E$ is the total energy. 

The collision operator $\mQ[f]$, including both $\mQ_{\text{Bol}}(f,f)$ and $\mQ_{\text{BGK}}(f)$, satisfies the conservation of mass, momentum, and energy:
\begin{equation}\label{conservation}
\int_{\mathbb{R}^d} \mQ[f] \,\rd v=\int_{\mathbb{R}^d} \mQ[f] v \,\rd v=\int_{\mathbb{R}^d} \mQ[f] |v|^2 \,\rd v=0,
\end{equation}
and the celebrated Boltzmann's H-theorem:
\begin{equation}
\int_{\mathbb{R}^d} \mQ[f] \ln f  \,\rd v \leq 0,
\end{equation}
with equality if and only if $ f $ reaches the equilibrium $f_{\text{eq}}$.

\section{Methodology}
\label{sec:method}

\subsection{Fourier Neural Operator (FNO)}

The Fourier Neural Operator (FNO) is designed to learn mappings between input and output fields directly in function space, aiming to approximate solution operators of partial differential equations (PDEs) \cite{LHLA2023, KLLABSA2023}. Its central idea is to perform learnable linear transformations on nonlocal operators in the frequency domain, combined with nonlinear activation functions. This approach enables approximation in infinite-dimensional spaces while maintaining robustness to discretization.

In the general set-up, given a training dataset $\mathcal{S}=\{(a^{(j)},u^{(j)})\}_{j=1}^n$, where the input functions $a \in \mathcal{A}: \mathcal{D} \mapsto \mathbb{R}$ and the output functions $u \in \mathcal{U}: \mathcal{D} \mapsto \mathbb{R}$ are defined on the same domain $\mathcal{D} \subset \mathbb{R}^d$, the goal is to approximate the true operator $\mathcal{G}: \mathcal{A} \mapsto \mathcal{U}$ by learning a parameterized operator $\G_\theta: \mathcal{A} \mapsto \mathcal{U}$, $\theta \in \Theta$. The optimization objective is to minimize the relative error:
\begin{equation} %
\min_{\theta \in \Theta} \  \sum_{(a,u) \in \mathcal{S}} \mathcal{L}(\G_\theta(a), u), 
\quad 
\end{equation}
where
\begin{equation}\label{loss1}
\mathcal{L}(\G_\theta(a),u) = \frac{\| \G_\theta(a) - u \|}{\| u \|}
\end{equation}
with $\|\cdot\|$ being a suitable norm.

To apply the FNO to solve the kinetic equation \eqref{kin}, we denote $\G_t: f_0(v) \mapsto f_t(v)$ as the evolution operator mapping the initial state $f(t=0,v) = f_0(v)$ to the state $f(t,v)=f_t(v)$ at time $t$, governed by the BGK and Boltzmann collision operator, respectively:
\begin{equation}
\left\{
\begin{aligned}
     &\partial_t f (t,v) = \mathcal{Q}[f](t,v), \ t>0, \ v \in \mathcal{D},\\
     &f(t=0,v) = f_0(v),
\end{aligned}
\right.
\end{equation}
then, in this particular case, we aim to approximate the evolution operator $\mathcal{G}_t$ by $\G_\theta$, i.e.,
\begin{equation}
  f_t(v) = \G_t[f_0](v) \approx \G_\theta[f_0](v).
\end{equation}
For example, if we take the norm in \eqref{loss1} as the $L^2$ Lebesgue norm, given a training dataset $\mathcal{S}=\{(f_0^{(j)},f_t^{(j)})\}_{j=1}^n$, the loss function reduces to 

\begin{equation}
\mathcal{L}(\G_\theta[f_0],f_t) = 
\frac{ \Big( \sum_{i=1}^N (\G_\theta[f_0](v_i) - f_t(v_i))^2 \Big)^{1/2}}
     {\Big( \sum_{i=1}^N f_t(v_i)^2 \Big)^{1/2}},
\end{equation}
where we suppose that the functions can be discretized on grids $\{v_i\}_{i=1}^N$.

Now, we are in a position to introduce how to specifically apply the standard FNO architecture  \cite{anandkumar2020neural} in the kinetic setting, which includes three modules: a lifting operator $L$, a stack of Fourier layers, and a projection operator $P$: 

\textit{Step I: Lifting.} First of all, the input function $f_0$ (concatenated with the velocity coordinate $v$ to capture velocity-dependent features) is lifted into a higher-dimensional feature space $F_0$ in the following form:
\begin{equation}
F_0(v) = L \big([\,v,\;f_0(v)\,]\big),
\end{equation}
where $L$ is typically a fully connected shallow network that expands the channel dimensions.   

\textit{Step II: Fourier layer.} In the Fourier layer, the update of the hidden representation $F_{t-1}$ can be summarized by the following sequence:
\[
F_{t-1}(v) \xrightarrow{\mathcal{F}} 
\hat{F}_{t-1}(k) \xrightarrow{\text{frequency update}} 
\tilde{F}_t(k) \xrightarrow{\mathcal{F}^{-1}} 
 \mathcal{F}^{-1}[\tilde{F}_t]\xrightarrow{\text{nonlinear}} 
F_t(v)
\]

First of all,
\begin{equation}
    \hat{F}_{t-1}(k) : = \mathcal{F}(F_{t-1}(v)),
\end{equation}
where $F_{t-1}(v)$ denotes the hidden features in the physical (velocity) domain at $t-1$, $\mathcal{F}$ denotes the discrete Fourier transform, and $\hat{F}_{t-1}(k)$ is the corresponding representation in the Fourier (frequency) domain at $t-1$. 

Then, $\tilde{F}_t(k)$ is the extended spectrum after the learnable Fourier transformation:
\begin{equation}
\tilde F_t(k):=
\left \{
\begin{aligned}
R_t(k)\,\hat F_{t-1}(k), & \qquad  k \in \mathbb{Z}_k,\\[2pt]
0, & \qquad \text{otherwise},
\end{aligned}
\right.
\end{equation}
where $\mathcal{F}^{-1}$ denotes the discrete inverse Fourier transform, $\mathbb{Z}_k$ is the set of retained low-frequency modes ($\|k\|_\infty \leq k_{\max}$), and $R_t(k) \in \mathbb{C}^{N_c \times N_c}$ are learnable matrices that mix channels in the frequency domain.

Next, $F_t(v)$ can be updated in the following form:
here each Fourier layer can be interpreted as approximating a kernel integral operator,
\begin{equation}\label{Ftv}
    F_{t}(v) = \sigma\Big( W_{t} F_{t-1}(v) + \mathcal{K}_\phi[ F_{t-1}](v) \Big), 
\end{equation}
where, $W_t$ is a linear transformation, $\sigma $ is a nonlinear activation function, 
\begin{equation}
    \mathcal{K}_\phi[F](v) := \int \kappa_\phi(v,v_*) F(v_*) \,\rd v_*.
\end{equation}
When $\kappa_\phi(v,v_*) = \kappa_\phi(v-v_*)$ is parametrized kernel with translation-invariant property, the operator $\mathcal{K}_\phi$ admits an efficient way to compute the convolution by the Fourier transform:
\begin{equation}\label{KF}
    \mathcal{K}_\phi[F](v):= (\kappa_\phi*F)(v) = \mathcal{F}^{-1} \Big( \tilde{F}(\kappa_\phi) \cdot \tilde{F}(v) \Big)(v),
\end{equation}
which illustrates that the Fourier layer of FNO is exactly a learnable approximation of such nonlocal collision operators.

\textit{Step III: Projection.} Finally, the hidden state is projected back to the output space:
\[
f_t(v)=P \big(F_t\big),
\]
where, $P$ denotes a projection operator that maps the high-dimensional latent representation $F_t$ back to the solution $f_t$ in physical space.

The three steps above describe the basic architecture of a FNO with one single Fourier layer. In practical implementations, multiple Fourier layers are typically stacked. Specifically, Step II is repeated for multiple layers, with the output of each Fourier layer serving as the input to the subsequent layer.

This architecture leverages Fourier domain multiplications to efficiently approximate nonlocal integral operators. By constraining the learnable parameters to low-frequency modes, it achieves discretization invariance and generalizes effectively across resolutions. The explicit inclusion of the coordinate variable $v$ enables the modeling of velocity-dependent effects, such as source terms and boundary conditions. In addition, the FNO weights are defined on the Fourier basis, which is globally supported on the continuous domain and independent of the sampling grid. Consequently, a model trained at one resolution can be directly evaluated at another, potentially enabling zero-shot super-resolution.

\subsection{FNO with conservative constraints (C-FNO)}

Recall that, in the standard FNO framework, the mean squared error (MSE) is used in the loss function, 
\begin{equation}
  \mathcal{L}_{\text{MSE}} := \frac{1}{n}\sum_{j=1}^n \| f_\theta^{(j)} - f_{\text{ref}}^{(j)} \|_{L^2}^2,
\end{equation}
where $\{f_\theta^{(j)}\}_{j=1}^n = \{\mathcal{G}_{\theta}^{(j)}[f_0]\}_{j=1}^n$ are the solutions obtained from the parametrized neural operator $\G_{\theta}$ at time $t$, and $\{f_{\text{ref}}^{(j)}\}_{j=1}^n =\{ f_{t}^{(j)}\}_{j=1}^n$ are the reference solutions at time $t$. 

However, both the Boltzmann and BGK collision operators satisfy the conservation of mass, momentum, and energy as in \eqref{conservation}.
It is natural to incorporate conservation penalties as soft constraints to achieve better conservative properties:
\begin{equation}\label{C-Loss}
  \mathcal{L}_{\text{Cons}}
  := \mathcal{L}_{\text{MSE}}
  + \lambda_m \mathcal{C}_m
  + \lambda_p \mathcal{C}_p
  + \lambda_e \mathcal{C}_e,
\end{equation}
with 
\begin{align*}
  \mathcal{C}_m &=
  \frac{1}{n}\sum_{j=1}^n \Big| \int_{\R^d} f_\theta^{(j)} \,\rd v - \int_{\R^d} f_{0}^{(j)} \,\rd v \Big|^2,\\[4pt]
  \mathcal{C}_p &=
  \frac{1}{n}\sum_{j=1}^n \Big\| \int_{\R^d} v\, f_\theta^{(j)} \,\rd v - \int_{\R^d} v\, f_{0}^{(j)} \,\rd v \Big\|_2^2,\\[4pt]
  \mathcal{C}_e &=
  \frac{1}{n}\sum_{j=1}^n \Big| \int_{\R^d} |v|^2 f_\theta^{(j)} \,\rd v - \int_{\R^d} |v|^2 f_{0}^{(j)} \,\rd v \Big|^2.
\end{align*}

With the newly proposed loss function with conservative constraints, the operator learning process of C-FNO is adapted to seek the minimizer of the following optimization problem:
\begin{equation}
\label{opt_cons}
  \min_{\theta \in \Theta}\;
  \mathbb{E}_{f_0}\!\left[
    \mathcal{L}_{\text{MSE}}
    + \lambda_m \mathcal{C}_m
    + \lambda_p \mathcal{C}_p
    + \lambda_e \mathcal{C}_e
  \right],
\end{equation}
where \(\mathcal{C}_m,\mathcal{C}_p,\mathcal{C}_e\) are the conservative penalties defined above. This refers to the \textit{point-to-point} learning framework throughout this paper, where the aim is to learn the evolution operator $\mathcal{G}_t$, i.e., the pair between the initial distribution function $f_0$ and the distribution function $f_t$ at time $t$. Observing from the following numerical experiments, these soft constraints turn out to promote physical consistency and enhance stability in long-term evolution.

In addition to learning the evolution operator from point to point, we also study the \textit{sequence-to-sequence} neural operator learning process, which is formulated in an iterative architecture, $\{f_0(v) \mapsto f_1(v) \mapsto \cdots \mapsto f_t(v) \mapsto \cdots \mapsto f_T(v) \}$,  where $\{f_t(v)\}$, defined for $t=1, 2, \cdots, T$, is the sequence of distribution functions in $v$. In this case, the loss function $\mathcal{L}_{\text{MSE}}$ in \eqref{opt_cons} is replaced with 
\begin{equation}
  \mathcal{L}_{\text{SS}} := \frac{1}{n} \sum_{j=1}^n \sum_{t=1}^T \| f_{t,\theta}^{(j)} - f_{t,\text{ref}}^{(j)} \|_{L^2}^2,
\end{equation}
where each sequence of $\{f_{1,\theta}^{(j)} \mapsto \cdots \mapsto f_{t,\theta}^{(j)} \mapsto \cdots \mapsto f_{T,\theta}^{(j)} \}$ are the solutions obtained from the parametrized neural operator, and
$\{f_{1,\text{ref}}^{(j)} \mapsto \cdots \mapsto f_{t,\text{ref}}^{(j)} \mapsto \cdots \mapsto f_{T,\text{ref}}^{(j)} \}$ are the reference data generated by the initial condition $\{f_0^{(j)}(v)\}_{j=1}^n$.

\subsection{Data generation and training process}

In this subsection, we present how the training dataset is prepared for both the BGK model and Boltzmann equation, which is crucial in the operator learning framework.

\subsubsection{Data generation and training for BGK model}
\label{subsubsec:data-BGK}

We consider the homogeneous BGK model \eqref{BGK} in both one-dimensional velocity and three-dimensional velocity space with the relaxation time $\text{Kn} = \tau=0.1$, and then $\mathcal{M}[f]$ is the Maxwellian distribution determined by the macroscopic quantities \eqref{rho-u-t} of $f$:
\[
\mathcal{M}[f](t,v)
=
\frac{\rho(t)}{(2\pi T(t))^{d/2}}
\exp\left(
-\frac{|v-u(t)|^2}{2T(t)}
\right).
\]
To generate the training data, we discretize the velocity domain \(v\in[-3,3]\) using 64 uniformly spaced grid points for one-dimensional case and using 32 uniformly spaced grid points for three-dimensional case. Then we advance the semi-discrete BGK equation with the forward Euler method using \(\Delta t=0.01\) for 100 time steps.
The macroscopic quantities are evaluated by the  standard trapezoidal rule:
\[
\rho^n=
\begin{cases}
\displaystyle
\sum_{i=1}^{N_v} w_i f_i^n,
& d=1,\\[3mm]
\displaystyle
\sum_{i=1}^{N_v}
\sum_{j=1}^{N_v}
\sum_{k=1}^{N_v}
w_i w_j w_k f_{ijk}^n,
& d=3,
\end{cases}
\]
\[
u^n=
\begin{cases}
\displaystyle
\frac{1}{\rho^n}
\sum_{i=1}^{N_v}
w_i v_i f_i^n,
& d=1,\\[3mm]
\displaystyle
\frac{1}{\rho^n}
\sum_{i=1}^{N_v}
\sum_{j=1}^{N_v}
\sum_{k=1}^{N_v}
w_i w_j w_k
v_{ijk} f_{ijk}^n,
& d=3,
\end{cases}
\]
\[
T^n=
\begin{cases}
\displaystyle
\frac{1}{\rho^n}
\sum_{i=1}^{N_v}
w_i (v_i-u^n)^2 f_i^n,
& d=1,\\[3mm]
\displaystyle
\frac{1}{3\rho^n}
\sum_{i=1}^{N_v}
\sum_{j=1}^{N_v}
\sum_{k=1}^{N_v}
w_i w_j w_k
\left|v_{ijk}-u^n\right|^2
f_{ijk}^n,
& d=3,
\end{cases}
\]
where, for $d=1$, the trapezoidal quadrature weights are defined as
\[
w_i=
\begin{cases}
\displaystyle
\frac{1}{2}\Delta v,
& i=1 \text{ or } i=N_v,\\[2mm]
\Delta v,
& \text{otherwise}.
\end{cases}
\]
and for $d=3$, the quadrature weight is
\[
w_{ijk}=w_iw_jw_k,
\]
with
\[
v_{ijk}
=
\begin{pmatrix}
v_i\\
v_j\\
v_k
\end{pmatrix}.
\]
Then, the Maxwellian distribution at $n$-th time step is assembled as 
\[
\mathcal{M}_i^n[f] \;=\; \frac{\rho^n}{(2\pi T^n)^{d/2}} \exp\Big( -\frac{|v_i-u^n|^2}{2 T^n} \Big),
\]
and training data is iterated as follows:
\begin{equation}
f_i^{n+1} \;=\; f_i^n + \Delta t\, \frac{\mathcal{M}_i^n - f_i^n}{\tau}.
\end{equation}
Note that for the parameter regime used in the BGK experiments ($\tau=0.1$ and $\Delta t=0.01$), a separate comparison with fourth order Runge–Kutta method (RK4) on the same velocity grid produced no material change in the reported solution profiles, supporting the use of forward Euler for reference data generation.

Following the algorithm above, we generate 100 independent trajectories (sequence of solutions) for both 1V and 3V cases respectively.
For each trajectory, we independently sample the initial macroscopic quantities
\[
\rho_0 \sim \mathcal{U}(0.8,\,1.2), \qquad
    u_0 \equiv 0, \qquad
    T_0 \sim \mathcal{U}(0.6,\,1.4),
\]
such that the initial distributions are constructed in the following diverse manner:
\begin{equation}\label{initial-BGK}
    f_{i}^{0,(j)} \;=\; 
\frac{\rho^{0,(j)}}{(2\pi T^{0,(j)})^{d/2}}
\exp\!\Big(-\tfrac{|v_i-u^{0,(j)}|^{2}}{2T^{0,(j)}}\Big)
\;(1+\; \varepsilon\,g(v_{i})^{(j)}),
\end{equation}
\begin{equation}
        g(v) = \sum_{k \in \{3,5,7,9,11,13\}} 
    a_k \cos\!\left(\frac{k\pi v}{L} + \phi_k\right), 
    v\in [-L,L]^d
\end{equation}
where $\varepsilon$ is the perturbation amplitude — larger $\varepsilon$ corresponding to stronger perturbations. \(\phi_k \sim \mathcal{U}(0,2\pi)\) are random phases and 
\(a_k \sim \mathcal{U}(0.4,1.0)\) with random signs \(\pm1\).
And the superscript $(j)$ indexes the trajectory ($j=1,\dots,100$) while $i$ indexes the velocity grid point.

\subsubsection{Data generation and training for Boltzmann equation} 
\label{subsec:data-Boltzmann}

For the Boltzmann equation, generating training data is particularly challenging, as solving the equation itself is a long-standing problem due to its high dimensionality, strong nonlinearity, and intrinsic nonlocality. To address this, we employ the fast Fourier spectral solver introduced in \cite{HQ2020} to obtain numerical solutions for use as training data. In applying the Fourier–Galerkin spectral method, we approximate the Boltzmann equation in the truncated computational domain $\mathcal{D}_L=[-L,L]^d$:
\begin{equation} \label{ABE}
\left\{
\begin{aligned}
&\partial_{t} f = \mathcal{Q}^{R}(f,f), \quad t>0, \ v\in \mathcal{D}_L,\\
& f(t= 0,v)= f_{0}(v), 
\end{aligned}
\right.
\end{equation}
where $f_{0}$ is a non-negative initial condition in $\mathcal{D}_L$. Here, $\mathcal{Q}^{R}$ is the truncated collision operator defined by
\begin{equation}\label{QR}
\begin{split}
\mathcal{Q}^{R}(g,f)(v)&=\int_{\mathcal{B}_R}\int_{\Sd^{d-1}}\Phi(|q|)b(\sigma\cdot \hat{q})\left[g(v_*')f(v')-g(v-q)f(v)\right]\, \rd{\sigma}\, \rd{q}\\
&=\int_{\mathbb{R}^d} \int_{\Sd^{d-1}}\mathbf{1}_{|q|\leq R}\Phi(|q|)b(\sigma\cdot \hat{q})\left[g(v_*')f(v')-g(v-q)f(v)\right]\,\rd{\sigma}\, \rd{q},
\end{split}
\end{equation}
where a change of variable $v_* \mapsto q=v-v_*$ is applied in \eqref{Bol}, and the new variable $q$ is truncated to a ball $\mathcal{B}_R$ centered at the origin with radius $R$. Here, we can write $q=|q|\hat{q}$ with $|q|$ being the magnitude and $\hat{q}$ being the direction over $\mathbb{S}^{d-1}$. Accordingly, 
\begin{equation}
v'=v-\frac{q-|q|\sigma}{2}, \quad v_*'=v-\frac{q+|q|\sigma}{2}.
\end{equation}
In practice, $L$ and $R$ are often chosen by an anti-aliasing argument \cite{PR00}: assume that $\text{Supp} (f_0(v))\subset \mathcal{B}_S$, then one can take
\begin{equation} \label{RL1}
R=2S, \quad L\geq \frac{3+\sqrt{2}}{2}S.
\end{equation}

Given $N\geq0$, we then seek a truncated Fourier series expansion of $f$ as
\begin{equation}
f(t,v)\approx f_N(t,v)=\sum\limits_{|k| = -N/2}^{N/2-1} \hat{f}_k(t) \e^{\im \frac{\pi}{L}k\cdot v} \in \mathbb{P}_N,
\end{equation}
where $k$ is a multi-index and 
\begin{equation}
\mathbb{P}_N=\text{span} \left\{ \e^{\im \frac{\pi}{L} k\cdot v}\Big| -N/2\leq |k| \leq N/2-1 \right\},
\end{equation}
equipped with inner product 
\begin{equation}
\langle f,g \rangle = \frac{1}{(2L)^{d}}\int_{\mathcal{D}_L} f \bar{g}\, \rd v.
\end{equation}
Substituting $f_N$ into (\ref{ABE}) and applying the Galerkin projection into the space $\mathbb{P}_N$ yields
\begin{equation} \label{PFS}
\left\{
\begin{aligned}
&\partial_{t} f_N(t,v) = \mP_N \mathcal{Q}^{R}(f_N,f_N)(t,v), \quad t>0, \ v\in \mathcal{D}_L,\\
& f_N(t=0,v)=f_{0,N}(v),
\end{aligned}
\right.
\end{equation}
where $\mP_N$ is the projection operator: for any function $g$,
\begin{equation}\label{proj}
\mP_N g=\sum_{|k|=-N/2}^{N/2-1} \hat{g}_k \e^{\im \frac{\pi}{L}k\cdot v}, \quad \hat{g}_k=\langle g, \e^{\im \frac{\pi}{L}k\cdot v}\rangle,
\end{equation}
and $f_{0,N} \in \mathbb{P}_N$ is a reasonable approximation to $f^0$, serving as the initial condition to the numerical system.

For each Fourier mode of \eqref{PFS}, we obtain,
\begin{equation} \label{FS}
\left\{
\begin{aligned}
&\partial_{t} \hat{f}_k(t) = \hat{\mathcal{Q}}^{R}_k(t), \quad  -N/2 \leq |k| \leq N/2-1,\\
& \hat{f}_k(t=0)= \hat{f}_{0,k},
\end{aligned}
\right.
\end{equation}
with
\begin{equation}
\hat{\mathcal{Q}}_{k}^R:=\langle \mathcal{Q}^R(f_N,f_N), \e^{\im \frac{\pi}{L}k\cdot v}\rangle, \qquad \hat{f}_{0,k}:=\langle f_{0,N}, \e^{\im \frac{\pi}{L}k\cdot v}\rangle.
\end{equation}
Using the definition in \eqref{QR} and orthogonality of the Fourier basis, we can derive that
\begin{equation} \label{sum}
\hat{\mathcal{Q}}_{k}^R =\sum\limits_{\substack{|l|,|m|=-N/2\\l+m=k}}^{N/2-1} G(l,m)\hat{f}_l \hat{f}_m,
\end{equation}
where the weight $G$ is given by
\begin{equation}\label{GG}
\begin{aligned}
G(l,m) &= \int_{\mathcal{B}_{R}}\int_{\Sd^{d-1}}\Phi(|q|)b(\sigma\cdot \hat{q})\left[ \e^{-\im \frac{\pi}{2L}(l+m)\cdot q +\im \frac{\pi}{2L}|q|(l-m)\cdot \sigma} - \e^{-\im \frac{\pi}{L}m\cdot q} \right]\,\rd\sigma\,\rd q\\
&=  \int_{\mathcal{B}_{R}}\e^{-\im \frac{\pi}{L}m\cdot q}\left[\int_{\Sd^{d-1}}\Phi(|q|)b(\sigma\cdot \hat{q})(\e^{\im \frac{\pi}{2L}(l+m)\cdot (q-|q|\sigma)}-1)\, \rd\sigma\right] \,\rd q. 
\end{aligned}
\end{equation}
The second equality follows from exchanging the variables $\sigma \leftrightarrow \hat{q}$ in the gain part of $G(l,m)$. 
Setting \(k=0\) in \eqref{GG} implies \(l+m=0\), for which \(G(l,-l)=0\). Consequently, \(\widehat Q^{R}_{0}=0\), and the zeroth Fourier coefficient of the spectral solution is preserved up to floating-point roundoff. Under the normalization in \eqref{proj}, the discrete mass is \(\rho_h(t)=(2L)^d\widehat f_0(t)\). We therefore evaluate the mass of the spectral solution from its zeroth Fourier coefficient, or equivalently from the Fourier-consistent periodic grid sum.
In the direct Fourier spectral method, $G(l,m)$ is precomputed, as it is independent of the solution, and the summation \eqref{sum} is then evaluated directly during the online computation. 
However, in the fast algorithm \cite{GHHH17, HQ2020}, if a low-rank decomposition of $G(l,m)$ in \eqref{GG} can be obtained in the form
\begin{equation}\label{GG1}
G(l,m)\approx\sum_{|p|=1}^{N_{p}}\alpha^{p}(k)\beta^{p}(m), \quad k = l + m,
\end{equation}
where $\alpha^{p}(k) = \alpha^{p}_k$ and $\beta^{p}(m) = \beta^{p}_m$ are suitable functions to be determined and the number of terms $N_{p}$ is small, then \eqref{sum} becomes
\begin{equation}
\mQ_{k}^R\approx \sum_{|p|=1}^{N_{p}}\alpha^{p}_k \sum\limits_{\substack{|l|,|m|=-N/2\\ l+m=k}}^{N/2-1}\hat{f}_{l}\left(\beta^{p}_m \hat{f}_{m}\right),
\end{equation}
where the inner summation is a convolution between $f_{l}$ and $\beta^{p}(m)f_{m}$. As a result, the total computational cost of evaluating $\mQ^R_{k}$ (for all $k$) is reduced from $O(N^{2d})$ to $O(N_{p}N^{d}\log N)$ through the use of a few Fast Fourier Transforms (FFTs).

\begin{remark}
It is worth noting the connection between the fast Fourier spectral method and the Fourier Neural Operator (FNO) in solving the Boltzmann equation. In particular, if the forward Euler method with time step $\Delta t$ is applied, the homogeneous Boltzmann equation can be approximated as
\begin{equation} 
\begin{aligned} 
f_{t+1}(v) =& f_t(v) + \Delta t\mathcal{Q}(f_t,f_t)(v)\\[3pt] \approx & f_t(v) + \Delta t\mathcal{F}^{-1}[\hat{\mathcal{Q}}_{k}^R](v)\\[3pt]
=& f_t(v) + \Delta t \mathcal{F}^{-1} \Big[\sum\limits_{\substack{|l|,|m|=-N/2\\l+m=k}}^{N/2-1} G(l,m)\hat{f}_{t,l} \hat{f}_{t,m}\Big](v)\\[3pt]
\approx & f_t(v) + \Delta t \mathcal{F}^{-1} \Big[\sum_p\alpha^p_k \Big((\beta^p_m\hat{f}_{t,m})*\hat{f}_{t,l}\Big)\Big](v), \end{aligned} 
\end{equation}
where the second term plays the role analogous to \eqref{KF} in FNO. The main distinction is that in the spectral approximation the nonlinearity is realized through the quadratic collision operator, whereas in FNO it is replaced by kernel convolutions. The absence of quadratic nonlinearity in FNO is, however, expected to be compensated by the activation function \eqref{Ftv} during training.
\end{remark}

The training data are therefore generated using the fast Fourier spectral method for the
spatially homogeneous Boltzmann equation with Maxwell molecules (i.e., $\gamma = 0$ in \eqref{Bb}). 
The velocity domain is
\[
    D_L=[-L,L]^d,
    \qquad
    L=\frac{3+\sqrt{2}}{2}S,
    \qquad S=2.5,
\]
Each velocity direction is discretized by
$N=32$ cell-centered points,
\[
    v_{j_k}=-L+\left(j_k-\frac12\right)\Delta v,
    \qquad
    \Delta v=\frac{2L}{N},
    \qquad j_k=1,\ldots,N.
\]
The reference solutions are computed in double precision by the RK4 method with $\Delta t=0.1$.

The initial conditions consist of three families with 300 ($d=2$) and 50 ($d=3$) samples in total: Gaussian, sum of Gaussians, and perturbed-Gaussian distributions. This sampling strategy is motivated by \cite{LJLH2025} but with fewer samples.
\begin{itemize}
\item \textbf{Gaussian distributions.}

The first family of Gaussian distributions is defined by
\begin{equation}
    M_{i,d}^{\mathrm G}(v_j)
    =
    \mathcal N_{h,d}\left[
        \frac{1}{(2\pi\sigma_{i,d}^{2})^{d/2}}
        \exp\left(
            -\frac{|v-c_{i,d}|^2}{2\sigma_{i,d}^{2}}
        \right)
    \right](v_j),
    \label{eq:gaussian-training-d}
\end{equation}
where $\mathcal{N}_{h,d}$ is a discrete normalization factor for unit mass.
The components of \(c_{i,d}\) are sampled independently from
\(\mathcal U[-0.5,0.5]\), while
\begin{equation}
    \sigma_{i,d}\sim
    \begin{cases}
        \mathcal U[0.3,0.5], & d=2,\\
        \mathcal U[0.5,0.8], & d=3.
    \end{cases}
    \label{eq:gaussian-width-d}
\end{equation}

\item \textbf{Sum of Gaussian distributions.}

In addition to the standard sum of Gaussian distributions, the second family also includes the approximations by the shifted and perturbed Bobylev–Krook–Wu (BKW) solution of the form in $d$-dimension:

The \(d\)-dimensional BKW solution is written as
\begin{equation}\label{BKW}
\begin{aligned}
    f_{\mathrm{BKW}}^{(d)}(t,v)
    &=
    \frac{1}{2(2\pi K_d(t))^{d/2}}
    \exp\left(-\frac{|v|^2}{2K_d(t)}\right)
    \left[
        \frac{(d+2)K_d(t)-d}{K_d(t)}
        +\frac{1-K_d(t)}{K_d(t)^2}|v|^2
    \right],
\end{aligned}
\end{equation}
where the time convention depends on the dimension:
\begin{equation}\label{Kd}
    K_d(t)=
    \begin{cases}
        1-\dfrac12\exp(-t/8), & d=2,\\[6pt]
        1-\exp(-t/6), & d=3.
    \end{cases}
\end{equation}

For $d=2$, all 100 training inputs are exact time-shifted states,
\begin{equation}
    M_{i,2}^{\mathrm S}(v)
    =
    f_{\mathrm{BKW}}^{(2)}(\tau_{i,2},v),
    \qquad
    \tau_{i,2}
    =
    \frac{2(i-\frac12)}{100},
    \qquad i=1,\ldots,100.
    \label{eq:2V-bkw-sampling}
\end{equation}
Thus, \(\tau_{i,2}=0.01,0.03,\ldots,1.99\), while the exact state at
\(t=0\) is excluded from training and reserved for the benchmark.

For \(d=3\), the 17 inputs form a hybrid family from the following time-shifted and perturbed initial conditions
\begin{equation}
    M_{i,3}^{\mathrm S}(v_j)
    =
    \mathcal N_{h,3}\left[
      f_{\mathrm{BKW}}^{(3)}(\tau_{i,3},\cdot)
    \right](v_j),
    \qquad
    \tau_{i,3}=5.52+0.06(i-1),
    \quad i=1,\ldots,9.
    \label{eq:3V-exact-bkw-sampling}
\end{equation}
The remaining eight are bounded perturbations at
\(t_{\star}=5.5\):
\begin{equation}
    M_{i,3}^{\mathrm S}(v_j)
    =
    \mathcal N_{h,3}\left[
      f_{\mathrm{BKW}}^{(3)}(t_{\star},v_j)
      \bigl(1+\eta_i q_i(v_j)\bigr)
    \right],
    \qquad i=10,\ldots,17,
    \label{eq:3V-perturbed-bkw-sampling}
\end{equation}
where
\[
    \{\eta_i\}_{i=10}^{17}
    =
    \{0.02,0.03,0.04,0.05,0.06,0.07,0.085,0.10\},
\]
and \(q_i\) is constructed from a random quadratic polynomial in $v_1/L,\ldots,v_d/L$. 

\item \textbf{Perturbed-Gaussian distributions.}

The third family can be expressed uniformly as
\begin{equation}
    M_{i,d}^{\mathrm P}(v_j)
    =
    \mathcal N_{h,d}\left[
      M_{i,d}^{\mathrm G}(v_j)
      \bigl(1+\varepsilon_{i,d}p_{i,d}(v_j)\bigr)
    \right],
    \label{eq:perturbed-gaussian-d}
\end{equation}
where \(p_{i,d}\) is constructed from a random quadratic expansion in the
normalized velocity coordinates.  More precisely,
\begin{equation}
    p_{i,d}(v)
    =
    a_{i,0}
    +\sum_{k=1}^{d}a_{i,k}\frac{v_k}{L}
    +\sum_{1\leq k\leq\ell\leq d}
      a_{i,k\ell}\frac{v_kv_\ell}{L^2}.
    \label{eq:random-quadratic-d}
\end{equation}
where
\[
d = 2: \quad a_{i,\alpha}\sim\mathcal N(0,1), \quad  \varepsilon_{i,d} \sim \mathcal U[0.02,0.10];\qquad 
d = 3: \quad a_{i,\alpha}\sim\mathcal U[-1,1], \quad \varepsilon_{i,d} = 0.2.
\]
\end{itemize}
Thus, the number of initial trajectories in the three families is
\begin{equation}
    \bigl(M_{\mathrm G}^{(d)}, M_{\mathrm S}^{(d)}, M_{\mathrm P}^{(d)}\bigr)
    =
    \begin{cases}
        (100,100,100), & d=2,\\
        (16,17,17), & d=3.
    \end{cases}
    \label{eq:family-counts-d}
\end{equation}

\section{Numerical experiments}
\label{sec:numerical}

In this section, we present the numerical tests of our proposed framework across the dimensions. 
Through extensive experiments, we compare the performance with various hyperparameters, e.g., the number of Fourier modes, hidden channels, layers of the network for the sequential prediction task in the case of BGK model, and make the benchmark tests and compare the conservation properties in the case of Boltzmann equation.

\subsection{One-dimensional case (1V) - BGK model}
\label{subsec:1V-BGK}

In this subsection, the numerical results for the BGK model in one dimension are presented, where we adopt a sequence-to-sequence training strategy, i.e., the neural operator is trained for the distribution functions at each time step of a sequence from the current state.

\subsubsection{Basic setup}
We first establish a baseline configuration by fixing a standard set of model parameters, which serves as a reference for subsequent experiments as follows:
the number of Fourier modes $N_{F} = 32$, the number of hidden channels $N_{C} = 64$, the number of layers $N_{L} = 4$. To systematically evaluate the sensitivity of the proposed framework, we then vary one parameter at a time while keeping all others fixed at their baseline values. This procedure enables a controlled assessment of the influence of each hyperparameter on the overall model performance.

The ground-truth solutions used for training and validation are generated following the procedure described in Section~\ref{subsubsec:data-BGK}, where the homogeneous BGK model is solved by an explicit Euler scheme with 64 velocity grid points and 100 time steps. Thus, we directly build on the same data generation setup previously introduced, ensuring consistency between training and testing datasets.

In this section, we focus on testing the proposed framework under two typical dynamical regimes --- \textit{near-equilibrium} and \textit{far-from-equilibrium} --- which differ in how far the initial state deviates from the local Maxwellian equilibrium distribution.

\begin{itemize}
    \item \textbf{Near-equilibrium case:} \\
    The system starts close to equilibrium, corresponding to small perturbations around the Maxwellian distribution. In this regime, we set the perturbation amplitude in the initial setup \eqref{initial-BGK} as a small value, e.g., $\varepsilon = 0.2$, which results in smooth and mild temporal evolution.

    \item \textbf{Far-from-equilibrium case:} \\
    The initial distribution deviates significantly from the Maxwellian, leading to strong non-equilibrium effects. This is modeled by setting the initial condition \eqref{initial-BGK} with a much larger perturbation amplitude, e.g., $\varepsilon = 3.0$, which produces highly nonlinear and rapidly varying dynamics.
\end{itemize}

These two regimes allow us to evaluate the performance of the model in learning both mild (\textit{near-equilibrium}) and strong (\textit{far-from-equilibrium}) non-equilibrium behaviors, providing a comprehensive understanding of its capability across different dynamical conditions.

The network was optimized by the Adam optimizer with an initial learning rate of $1\times 10^{-3}$, decayed by a factor of $0.5$ every $50$ epochs. The total number of training epochs was set to $200$, with a batch size of $32$.
For all BGK models, the experiments are conducted on a MacBook Pro equipped with an Apple M2 chip, featuring an 8-core CPU (4 performance cores and 4 efficiency cores) and 8 GB of unified memory. The system was running macOS with kernel version 8422.141.2 and Darwin version 7459.141.1.700.1. Despite the limited memory, the model training and evaluation processes are completed efficiently, demonstrating the computational feasibility of the proposed approach even on resource-constrained consumer hardware.

\subsubsection{Numerical experiment 1: Sensitivity of Fourier modes}

In practice, the training model \eqref{Ftv} cannot apply all Fourier modes, instead, it only learns and applies transformations to a fixed number of low-frequency modes. The number of Fourier modes $N_F$ determines how many frequency components the network keeps and learns from.
In this subsection, we investigate the impact of the number of Fourier modes $N_F$ on the performance of the network.\\

\begin{figure}[H]
    \centering
    \includegraphics[width=1.03\textwidth]{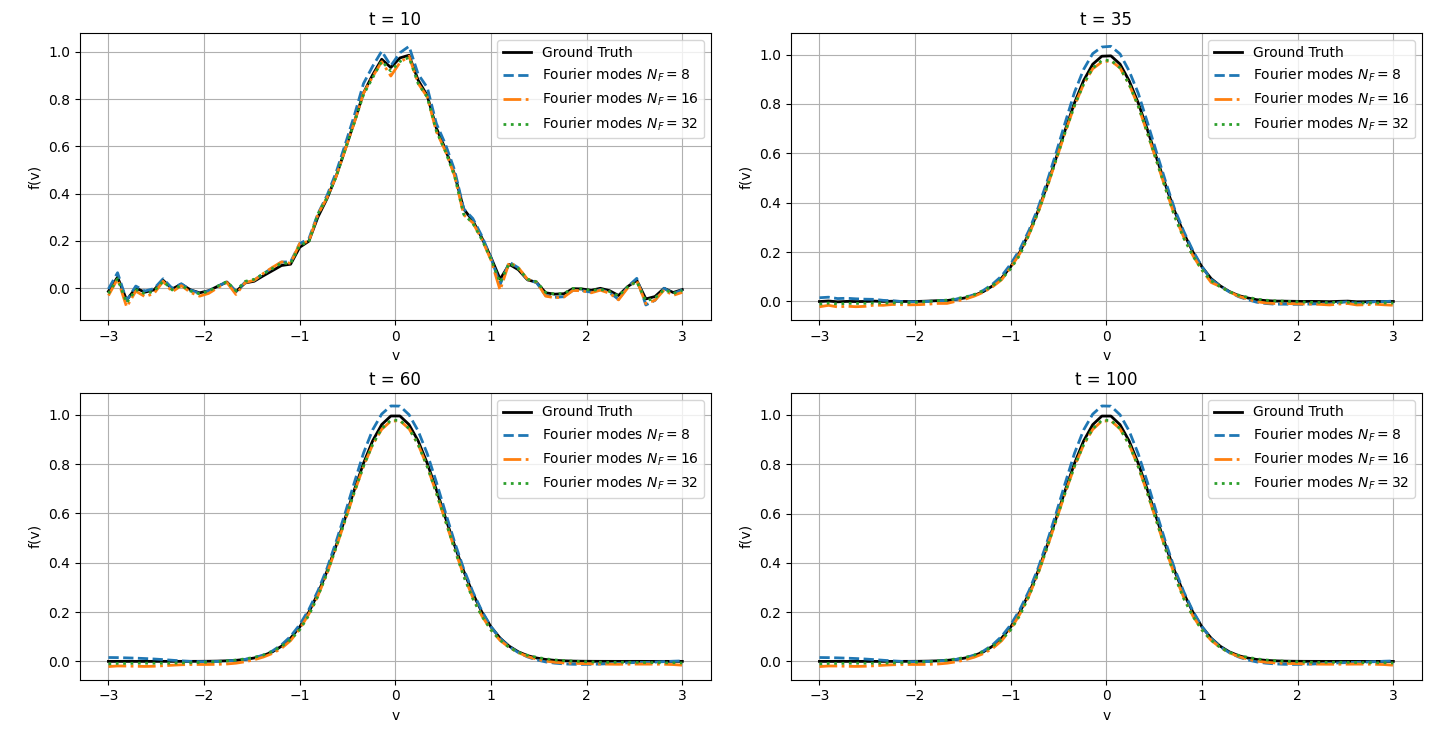}
    \caption{Evolution of $f(v)$ of the BGK model with different $N_F=8,16,32$ for near-equilibrium initial condition.}
    \label{fig:1V-near-sol}
\end{figure}

\steptitle{(a) Case 1.1: Near-equilibrium initial conditions for different $N_F$} 
In Figure~\ref{fig:1V-near-sol}, we examine the temporal evolution of the solutions obtained using different numbers of Fourier modes. The solutions are generated by the trained network at time steps $t = 10, 35, 60,$ and $100$ for a new initial condition close to equilibrium. All three configurations produce satisfactory results, with $N_F= 16$ and $N_F = 32$ yielding solutions that are nearly indistinguishable from the reference solution for all time. In contrast, the case with $N_F = 8$ exhibits slightly increased fluctuations during the early stages of evolution, although the overall accuracy remains acceptable.

\begin{figure}[H]
    \centering
    \includegraphics[width=0.7\textwidth]{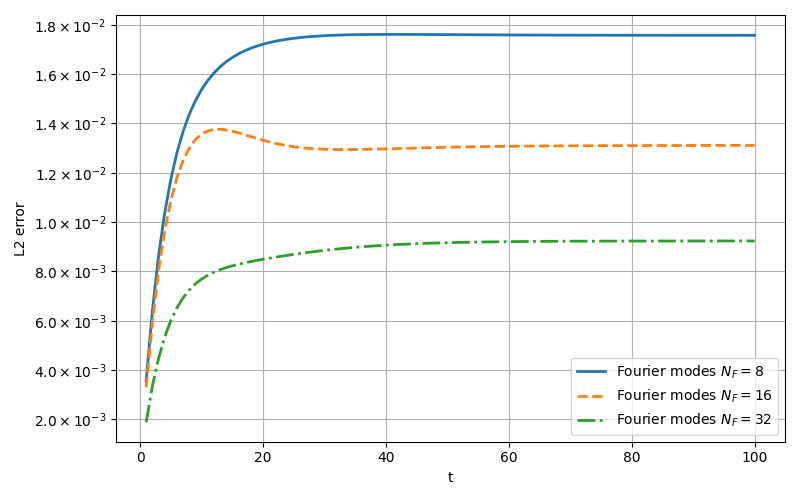}
    \caption{Evolution of $L^2$-error with $N_F=8,16,32$ for near-equilibrium initial condition.}
    \label{fig:1V-near-error-L2}
\end{figure}

\begin{figure}[H]
    \centering
    \includegraphics[width=0.7\textwidth]{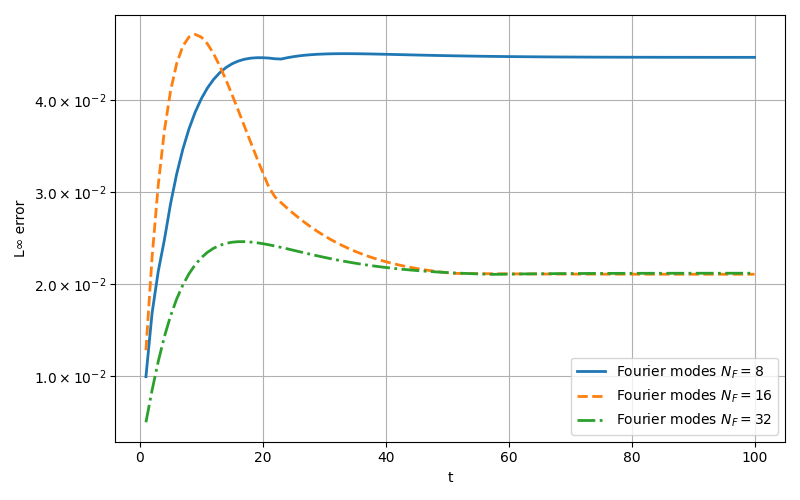}
    \caption{Evolution of $L^{\infty}$-error with $N_F=8,16,32$ for near-equilibrium initial condition.}
    \label{fig:1V-near-error-Linf}
\end{figure}

In Figure~\ref{fig:1V-near-error-L2} and Figure~\ref{fig:1V-near-error-Linf}, the temporal evolution of the error $\|f_{t,\theta}-f_{t,\text{Ref}}\|$ in both $L^2$ and $L^{\infty}$ norms are presented.

\steptitle{(b) Case 1.2: Far-from-equilibrium initial conditions for different $N_F$} 

\begin{figure}[H]
    \centering
    \includegraphics[width=1.03\textwidth]{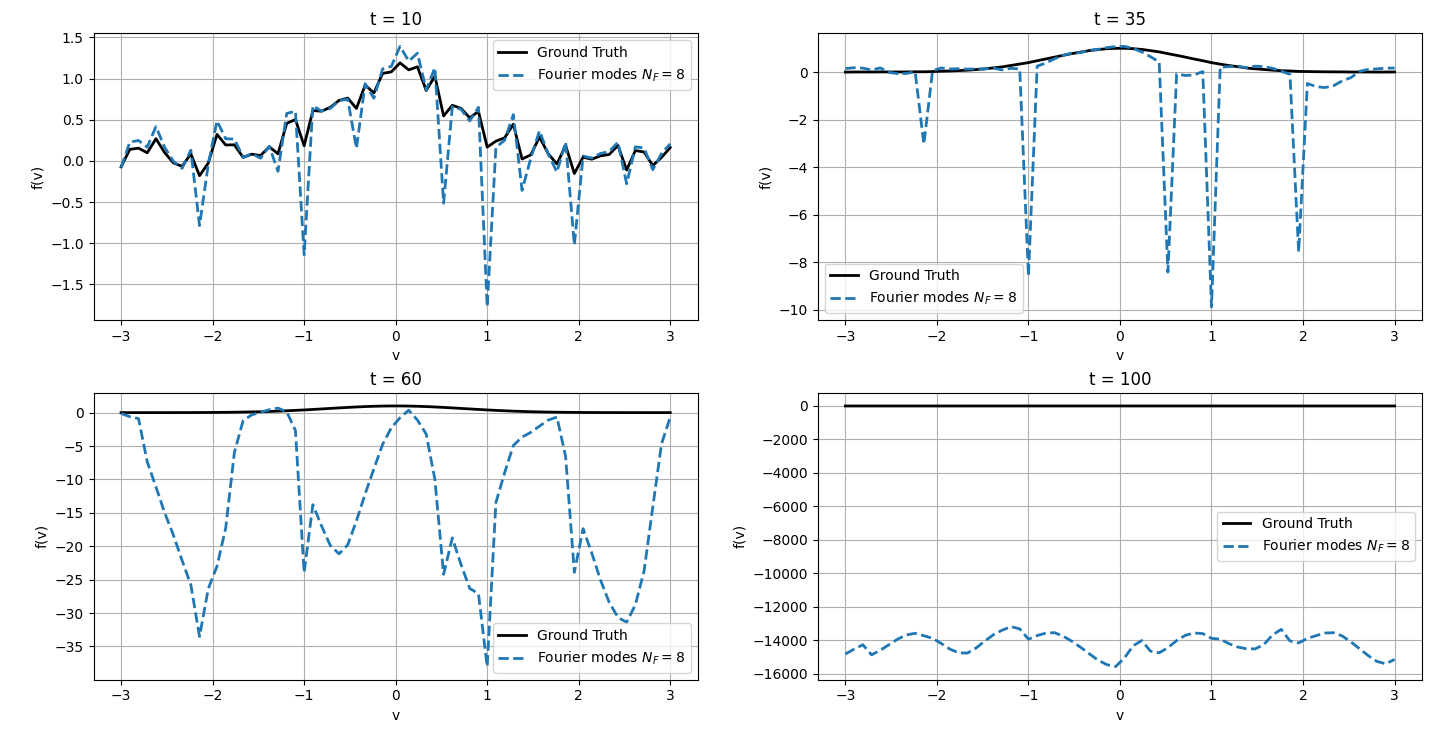}
    \caption{Evolution of $f(v)$ of the BGK model with $N_F=8$ for far-from-equilibrium initial condition.}
     \label{fig:1V-far-sol0}
\end{figure}

\begin{figure}[H]
    \centering
    \includegraphics[width=1.03\textwidth]{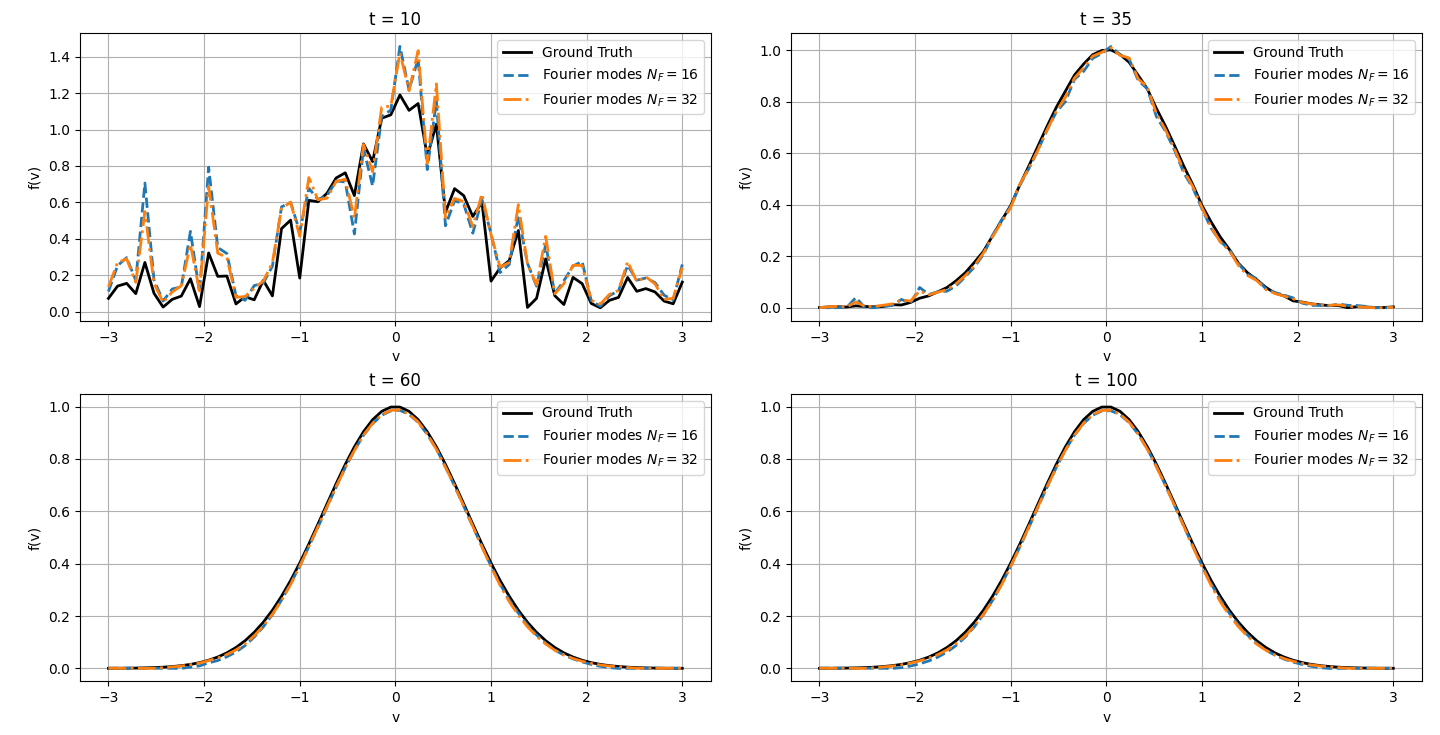}
    \caption{Evolution of $f(v)$ of the BGK model with different $N_F=16,32$ for far-from-equilibrium initial condition.}
    \label{fig:1V-far-sol1}
\end{figure}

In Figure~\ref{fig:1V-far-sol1}, we investigate the temporal evolution of the solutions obtained using different numbers of Fourier modes for far-from-equilibrium initial conditions. The configurations with $N_F = 16$ and $N_F = 32$ continue to yield stable and accurate results. In contrast, both insufficient and excessive numbers of modes can introduce undesirable oscillations in the predicted solutions, though for different underlying reasons. When the number of modes is too small (e.g., $N_F = 8$ in Figure~\ref{fig:1V-far-sol0}), the severe truncation in Fourier space causes the loss of significant high-frequency information, leading to noticeable deviations from the reference solution. Conversely, when the number of modes is excessively large, the network may overfit the training data, resulting in localized instabilities when evolving new initial conditions not encountered during training. These observations highlight the importance of selecting an appropriate truncation in the Fourier domain to balance the accuracy and generalization in the standard FNO framework.

\begin{figure}[H]
    \centering
    \includegraphics[width=0.7\textwidth]{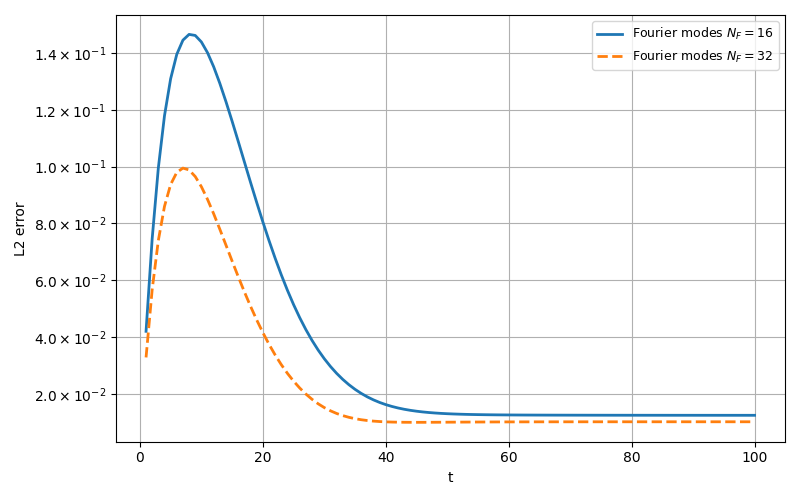}
    \caption{Evolution of $L^2$-error with $N_F=16,32$ for far-from-equilibrium initial condition.}
    \label{fig:1V-far-error-L2}
\end{figure}

\begin{figure}[H]
    \centering
    \includegraphics[width=0.7\textwidth]{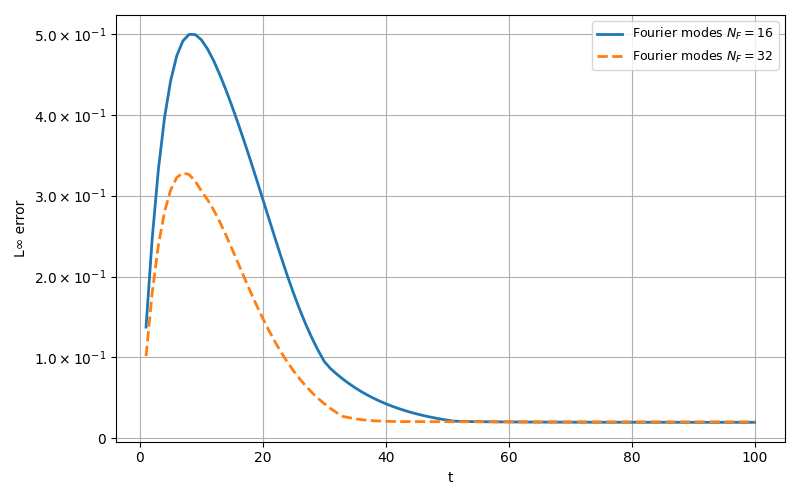}
    \caption{Evolution of $L^\infty$-error with $N_F=16,32$ for far-from-equilibrium initial condition.}
    \label{fig:1V-far-error-Linf}
\end{figure}

Furthermore, the temporal evolution of the error $\|f_{t,\theta}-f_{t,\text{Ref}}\|$ in both $L^2$ and $L^{\infty}$ norms are presented in Figure~\ref{fig:1V-far-error-L2} and Figure~\ref{fig:1V-far-error-Linf}.

\steptitle{(c) Computational cost for different $N_F$} 
In addition to evaluating the accuracy, we also assess the computational cost associated with the number of Fourier modes. The average training time per epoch for each configuration is summarized in Table~\ref{tab:runtime}. As anticipated, models with larger $N_F$ exhibit a higher computational overhead, primarily due to the increased cost of the Fourier transform and the larger number of trainable parameters involved in processing additional Fourier modes.
\begin{table}[H]
    \centering
    \begin{tabular}{ccc}
        \toprule
        \textbf{Modes} & \textbf{Training Parameters} (million) & \textbf{Training Time} (seconds) \\
        \midrule
        8   & 0.18&16.81 \\
        16  & 0.31&17.33 \\
        32  &0.57& 18.54 \\
        64 &1.1& 21.19 \\
        \bottomrule
    \end{tabular}
    \caption{Average training time per epoch for different numbers of Fourier modes $N_F = 8, 16, 32, 64$.}
    \label{tab:runtime}
\end{table}

\subsubsection{Numerical experiment 2: Sensitivity of hidden channels}

Recall that hidden channels refer to the number of feature dimensions (the ``width" of the feature space) used to represent the function after lifting from the input space, in this subsection, we compare the performance of the network in terms of the different hidden channel sizes $N_C$.\\

\steptitle{(a) Case 2.1: Near-equilibrium initial condition for different $N_C$}

\begin{figure}[H]
    \centering
    \includegraphics[width=01.03\textwidth]{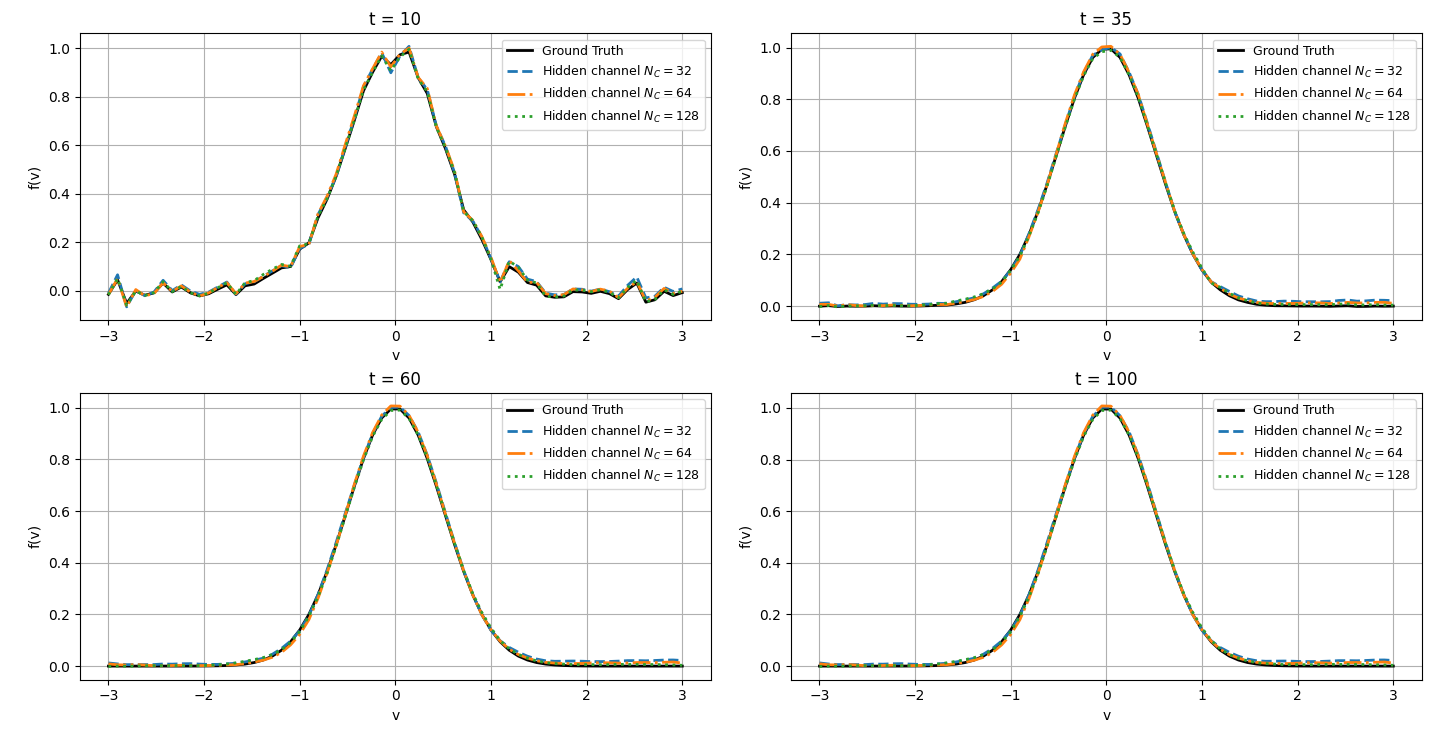}
    \caption{Evolution of $f(v)$ of the BGK model with different $N_C=32,64,128$ for near-equilibrium initial condition.}
    \label{fig:ch1}
\end{figure}

As shown in Figure~\ref{fig:ch1}, under the near-equilibrium scenario, all model variants—regardless of the number of hidden channels—perform very well. The predicted distributions align closely with the ground truth at every time step, and the discrepancy is barely noticeable. This indicates that for simple or smooth dynamics, increasing model capacity yields limited performance gain and may even be computationally inefficient.

\begin{figure}[H]
    \centering
    \includegraphics[width=0.7\textwidth]{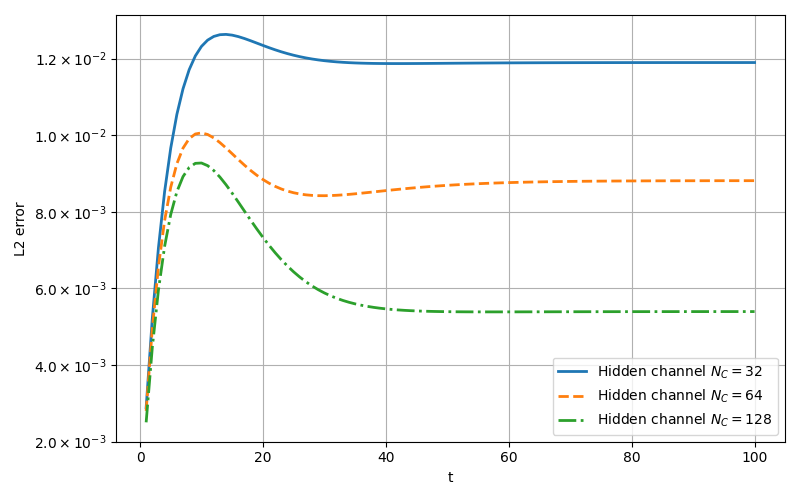}
    \caption{Evolution of $L^2$-error with $N_C=32,64,128$ for near-equilibrium initial condition.}
\end{figure}

\begin{figure}[H]
    \centering
    \includegraphics[width=0.7\textwidth]{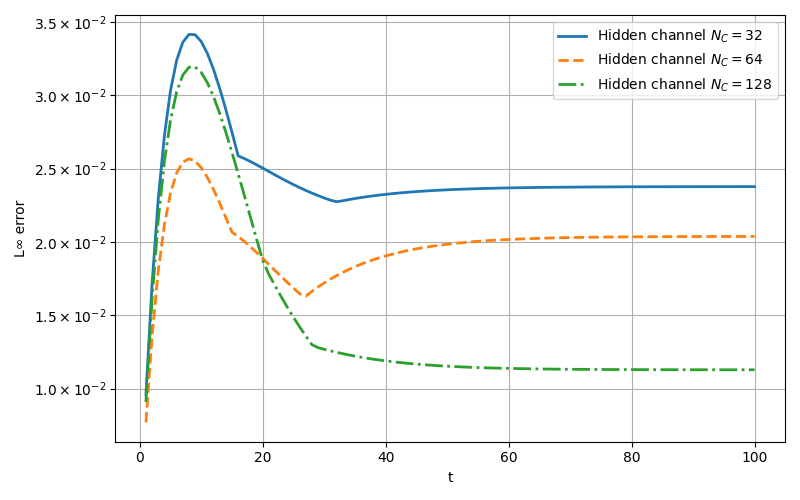}
    \caption{Evolution of $L^\infty$-error with $N_C=32,64,128$ for near-equilibrium initial condition.}
\end{figure}

\begin{figure}[H]
    \centering
    \includegraphics[width=1.03\textwidth]{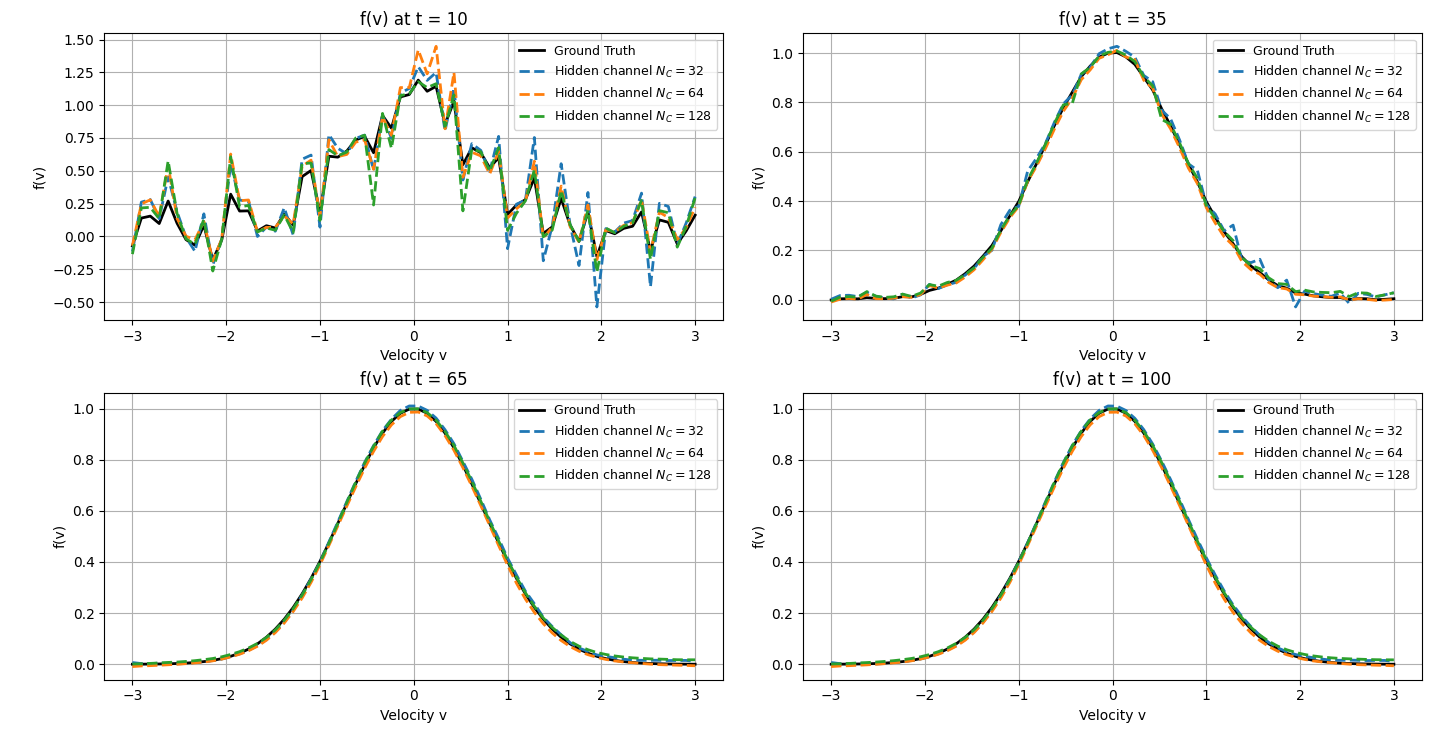}
    \caption{Evolution of $f(v)$ of the BGK model with different $N_C=32,64,128$ for far-from-equilibrium initial condition.}
    \label{fig:hidden2}
\end{figure}

\steptitle{(b) Case 2.2: Far-from-equilibrium initial condition for different $N_C$}
Figure~\ref{fig:hidden2} presents the temporal evolution of the velocity distribution $f(v)$ at selected time points ($t=10, 35, 65, 100$) for different hidden channel sizes. 
At the early stage ($t=10$), the solutions exhibit noticeable oscillations and deviations from the ground truth. 
Although the prediction errors are relatively large in this non-equilibrium regime, they remain within an acceptable range. 
As time progresses ($t=35, 65, 100$), the predicted distributions gradually converge toward the equilibrium, with the discrepancies among different hidden channel configurations becoming less pronounced. 
Overall, increasing the number of hidden channels enhances stability and accuracy during the transient phase, while all models achieve comparable performance in long-time regime.

\begin{figure}[H]
    \centering
    \includegraphics[width=0.7\textwidth]{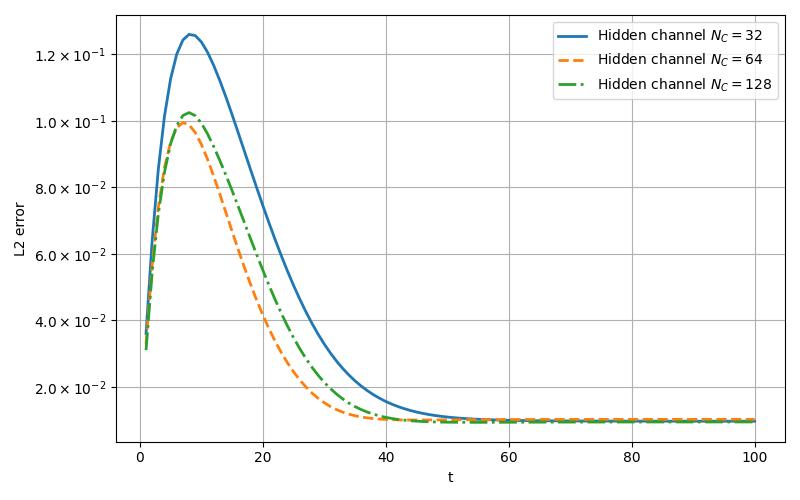}
    \caption{Evolution of $L^2$-error with $N_C=32,64,128$ for far-from-equilibrium initial condition.}
\end{figure}

\begin{figure}[H]
    \centering
    \includegraphics[width=0.7\textwidth]{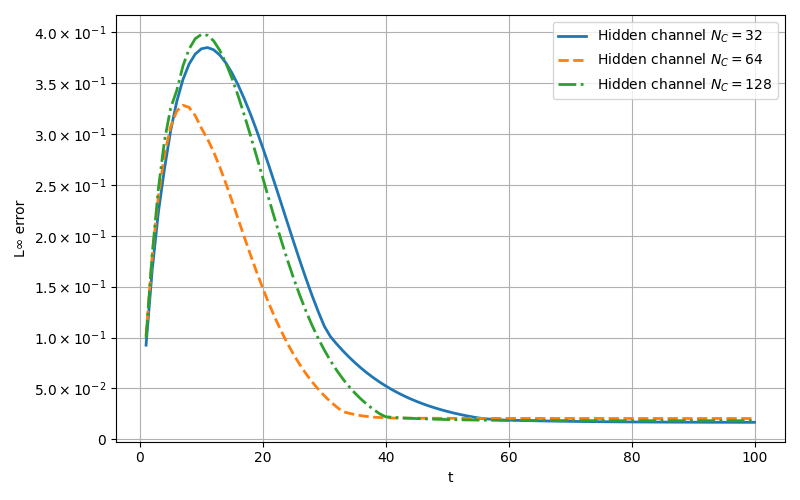}
    \caption{Evolution of $L^\infty$-error with $N_C=32,64,128$ for far-from-equilibrium initial condition.}
\end{figure}

\steptitle{(c) Computational cost for different $N_C$}  
As shown in Table~\ref{tab:runtime-2}, the computational cost grows significantly with the number of hidden channels $N_C$. Training of the model with $N_C = 32$ is completed in just 10.13 seconds, while the model with $N_C = 128$ requires 42.80 seconds. This sharp rise in training time reflects the increase in parameter count and the associated cost of spectral convolutions.

\begin{table}[H]
\centering
\begin{tabular}{ccc}
\toprule
\textbf{Hidden Channels $N_C$} & \textbf{Training Parameters} (million)& \textbf{Training Time} (seconds) \\
\midrule
32  &0.14& 10.13\\
64  &0.57& 18.54 \\
128 & 2.3&42.80\\
\bottomrule
\end{tabular}
\caption{Average training time per epoch for different numbers of hidden channels $N_C = 32,64,128$.}
\label{tab:runtime-2}
\end{table}

\subsubsection{Numerical experiment 3: Sensitivity of layers}

The number of layers $N_L$ represents the number of Fourier layers inserted in the neural operator learning framework (see \cite[Figure 2]{Li-ICLR2021}), which represents the depth of the architecture. 
In this subsection, we investigate how the number of layers $N_L$ affects performance.\\

\steptitle{(a) Case 3.1: Near-equilibrium initial condition for different $N_L$}

\begin{figure}[H]
    \centering
    \includegraphics[width=1.03\textwidth]{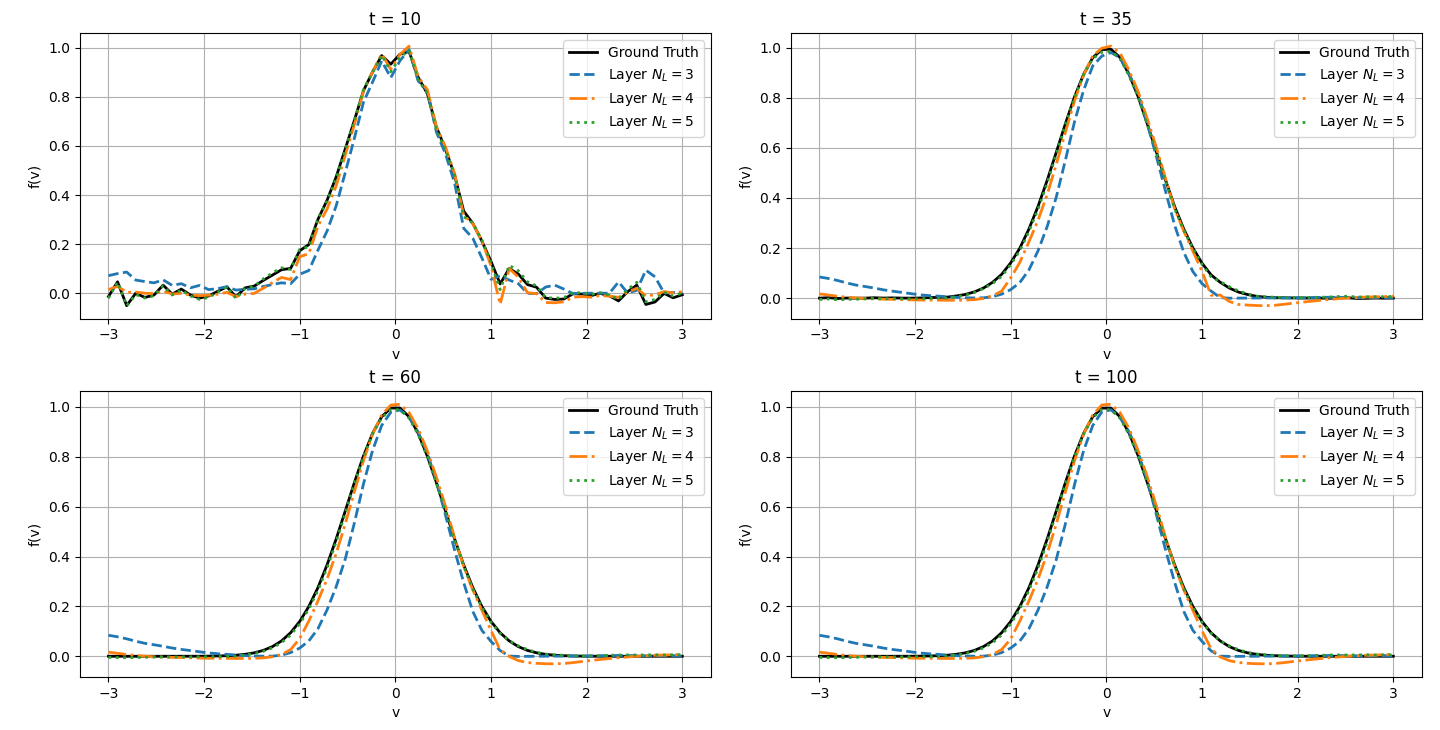}
    \caption{Evolution of $f(v)$ of the BGK model with different $N_L=3,4,5$ for near-equilibrium initial condition.}
    \label{layer1}
\end{figure}

As shown in Figure~\ref{layer1}, most models perform well under near-equilibrium conditions, with predicted profiles closely matching the ground truth across time steps. However, the model with only $N_L = 3$ layers clearly cannot perform well: it fails to capture the detailed structure of the distribution function, particularly at the early stage ($t=10$), and exhibits noticeable deviations at later times ($t=35$ and $t=60$). In contrast, models with 4 or 5 layers yield significantly more accurate and stable predictions, demonstrating excellent agreement with the reference solution.
These results suggest that although near-equilibrium dynamics are relatively simple to model, a minimum network depth of four layers is required to provide sufficient expressive capacity and prevent underfitting. Increasing the depth from $N_L = 4$ to $N_L = 5$ layers yields only marginal improvement, indicating decrease in returns. Hence, while deeper architectures are unnecessary in this regime, excessively shallow networks remain inadequate even for near-equilibrium scenarios.

\begin{figure}[H]
    \centering
    \includegraphics[width=0.7\textwidth]{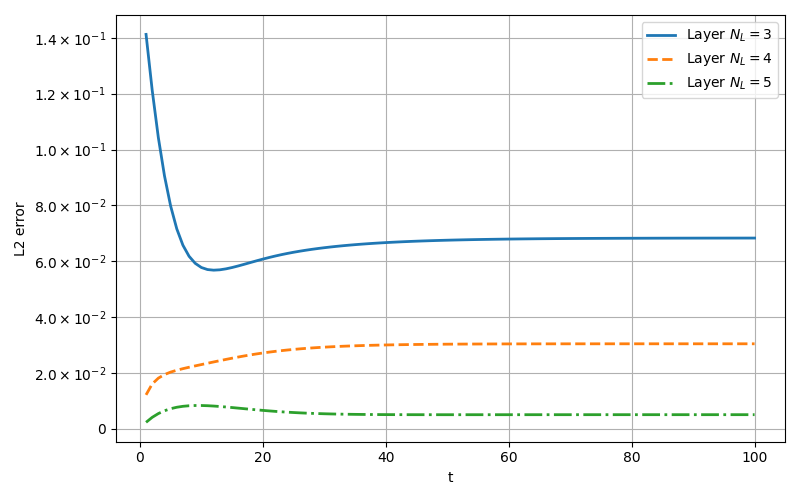}
    \caption{Evolution of $L^2$-error with $N_L=3,4,5$ for near-equilibrium initial condition.}
\end{figure}

\begin{figure}[H]
    \centering
    \includegraphics[width=0.7\textwidth]{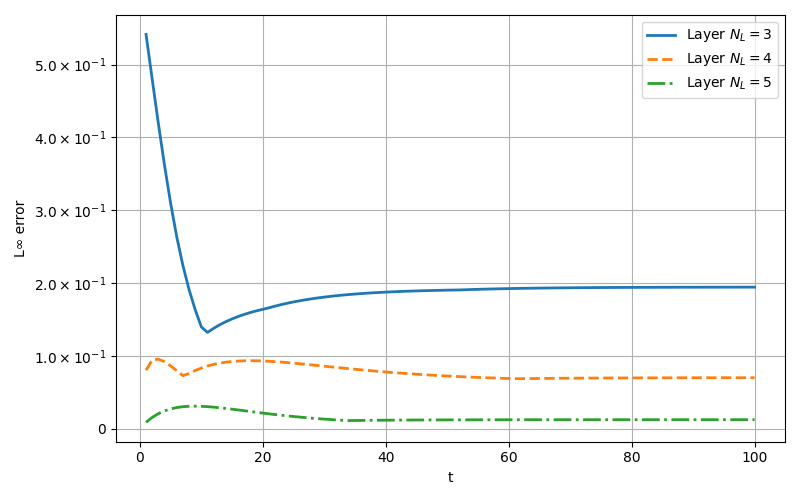}
    \caption{Evolution of $L^\infty$-error with $N_L=3,4,5$ for near-equilibrium initial condition.}
\end{figure}

\steptitle{(b) Case 3.2: Far-from-equilibrium initial condition for different $N_L$}

When the initial condition is far from the equilibrium state, the 3-layer model demonstrates clear deficiencies, particularly in early and mid time steps, producing significant oscillatory artifacts and underestimation in regions of sharp gradients. These issues are mostly mitigated when increasing the depth to $N_L = 4$ and further to $N_L = 5$. The 5-layer model consistently provides closer alignment to the ground truth, particularly at early time steps ($t=10$, $t=35$). It captures sharper transitions and avoids the smoothing effects observed in the 4-layer version. This suggests that increasing depth can enhance the ability of the learned operator to represent complex non-equilibrium dynamics without sacrificing stability.

\begin{figure}[H]
    \centering
    \includegraphics[width=1.03\textwidth]{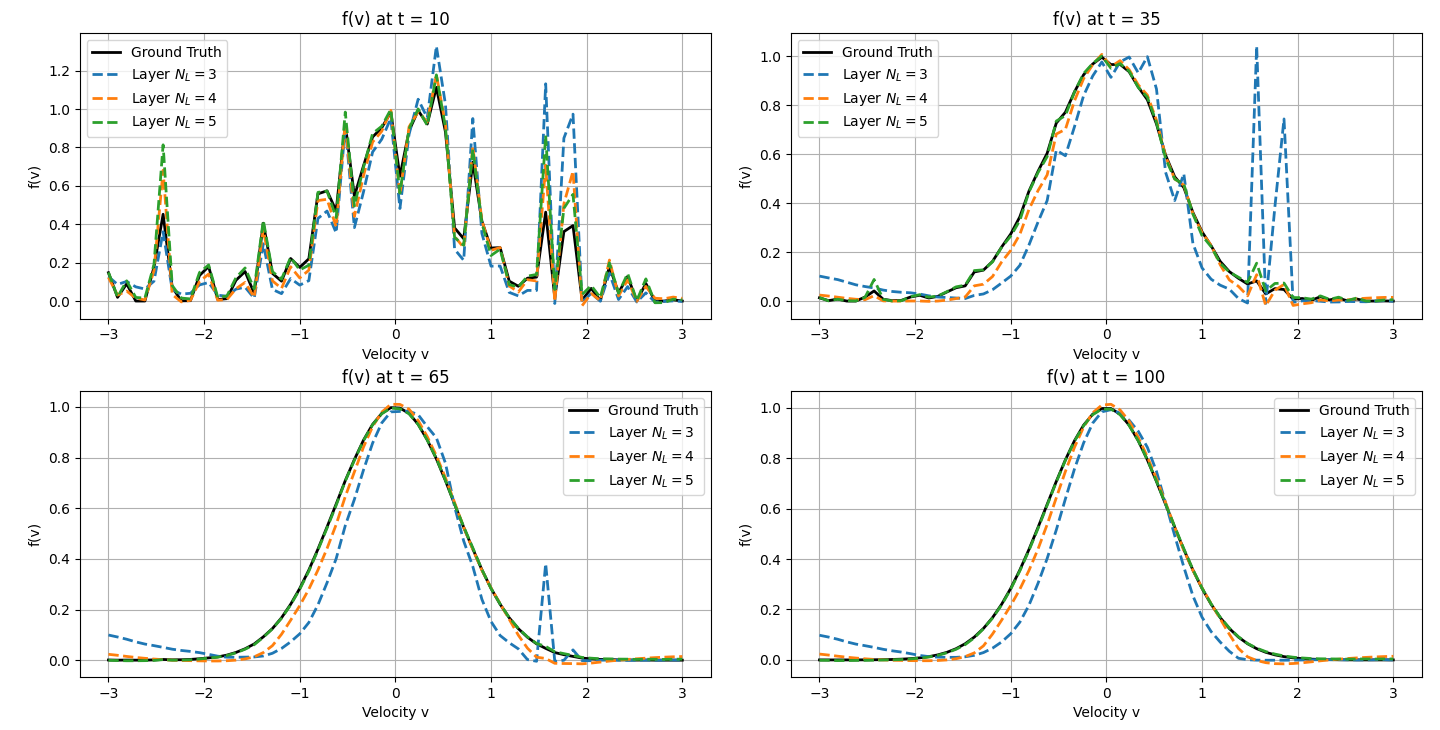}
    \caption{Evolution of $f(v)$ of the BGK model with different $N_L=3,4,5$ for far-from-equilibrium initial condition.}
    \label{layer2}
\end{figure}

\begin{figure}[H]
    \centering
    \includegraphics[width=0.7\textwidth]{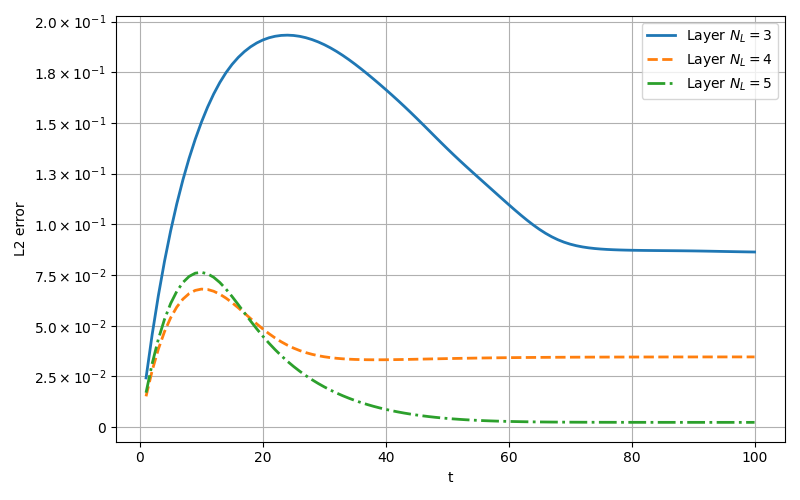}
    \caption{Evolution of $L^2$-error with $N_L=3,4,5$ for far-from-equilibrium initial condition.}
\end{figure}

\begin{figure}[H]
    \centering
    \includegraphics[width=0.7\textwidth]{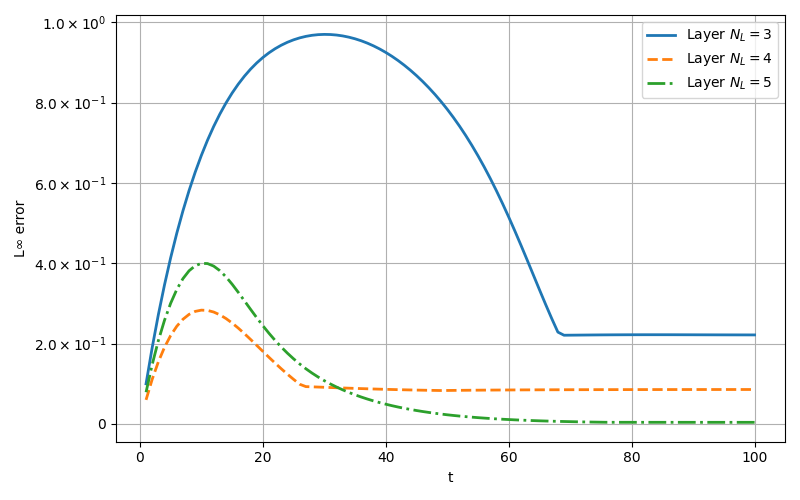}
    \caption{Evolution of $L^\infty$-error with $N_L=3,4,5$ for far-from-equilibrium initial condition.}
\end{figure}

\steptitle{(c) Computational cost for different $N_L$}  

To evaluate the computational cost associated with different depths, we measured the average training time per epoch for models with $N_L=3,4,5$. As shown in Table~\ref{tab:runtime_layers}, all configurations demonstrate comparable efficiency. These findings indicate that increasing depth does not necessarily lead to significantly higher training cost within this range.
\begin{table}[H] 
\centering 
\begin{tabular}{ccc} 
\toprule 
\textbf{Layers $N_L$} & \textbf{Training Parameters} (million) & \textbf{Training Time} (seconds) \\ 
\midrule 
3 & 0.43&15.10 \\
4 & 0.57&18.54 \\
5 & 0.72&21.96 \\
\bottomrule 
\end{tabular} 
\caption{Average training time per epoch for different numbers of layers $N_L = 3,4,5$.} 
\label{tab:runtime_layers} 
\end{table}

\subsubsection{Numerical experiment 4: Comparison with other operator learning framework}

\begin{figure}[H]
    \centering
    \includegraphics[width=1\textwidth]{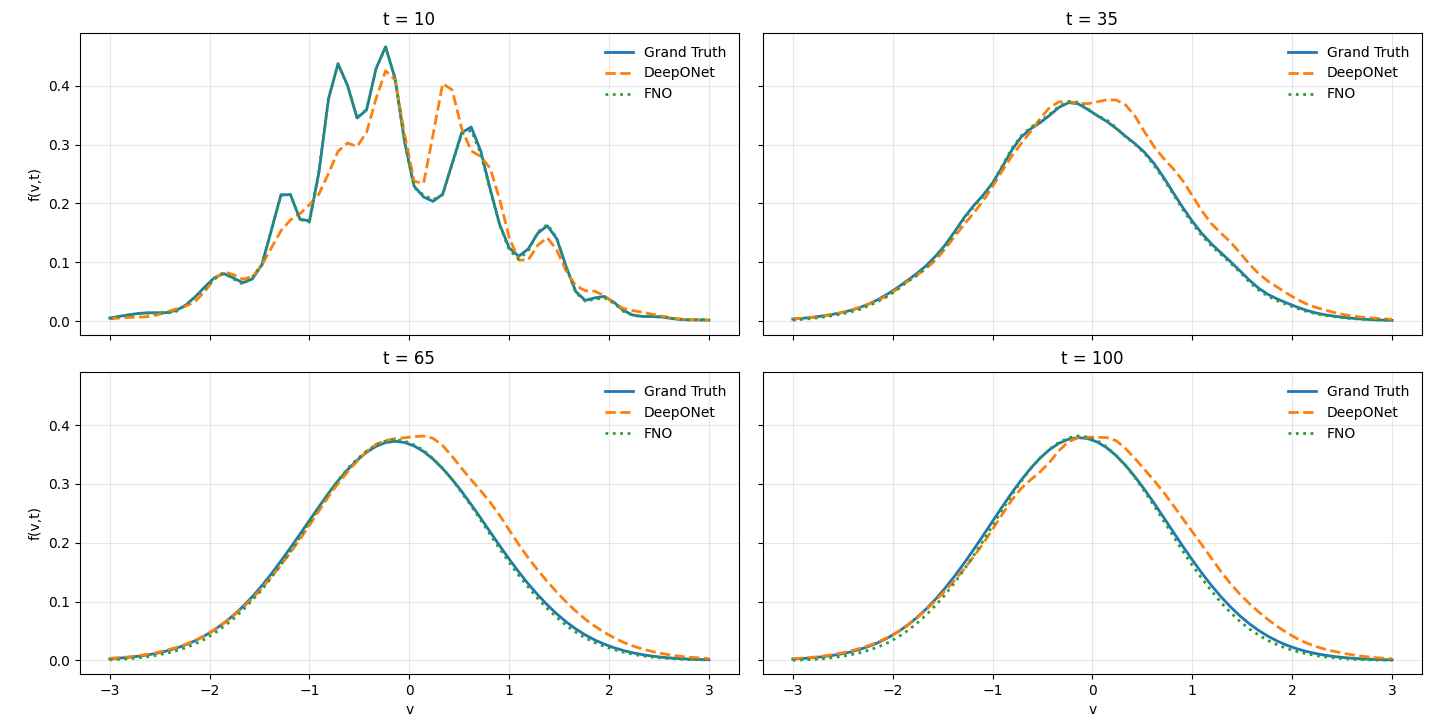}
    \caption{Evolution of $f(v)$ of the BGK model by FNO and DeepONet }
    \label{fig:DeepOnet}
\end{figure}

\begin{figure}[H]
    \centering
    \includegraphics[width=0.7\textwidth]{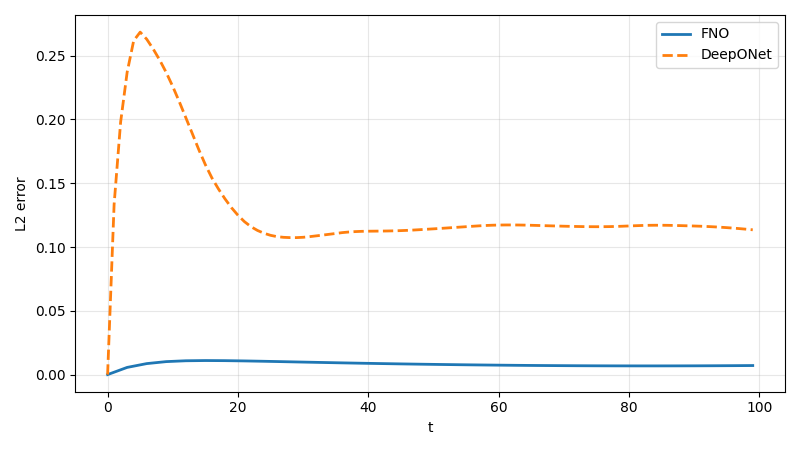}
    \caption{Evolution of $L^2$-error by FNO and DeepONet}
    \label{fig:DeepOnet_L2}
\end{figure}

In this experiment, we provide a comparative benchmark against DeepONet \cite{lu2021learning}, a widely used operator-learning framework. For the FNO model, we use the hyperparameters selected from the preceding sensitivity studies, namely \(N_F=32\), \(N_C=64\), and \(N_L=4\). For DeepONet model, the layer sizes for branch network and trunk network are $[64,512,512,512]$ and $[1,512,512,512]$, respectively, with ReLU activation and Glorot normal initialization. And FNO is trained by using the Adam optimizer, using learning rate \(10^{-3}\), batch size \(32\), and \(100\) training epochs. DeepONet is first trained using the Adam optimizer with a learning rate of
$10^{-3}$ for 200 steps and is subsequently refined using L-BFGS until
convergence. The complete training process consists of 50,633 optimization steps. Figure~\ref{fig:DeepOnet} compares the predicted velocity distributions with the reference solution. The DeepONet prediction exhibits visible deviations over the velocity domain, whereas the FNO prediction remains closer to the reference solution throughout the relaxation process. The corresponding time evolution of the \(L^2\)-error is shown in Figure~\ref{fig:DeepOnet_L2}, where the FNO error remains consistently smaller than that of DeepONet.


Table~\ref{tab:runtime-3} further reports the computational cost. Although DeepONet has a smaller per-step runtime, it requires substantially more optimization steps to reach a comparable level of convergence in our experiments. In contrast, FNO achieves better predictive accuracy with fewer training epochs, leading to a smaller total training time in this test. These results indicate that the Fourier-domain architecture is particularly effective for the present kinetic relaxation problem.

\begin{table}[H]
\centering
\begin{tabular}{cccc}
\toprule
\textbf{Method} & \textbf{Training Parameters} (million) 
& \textbf{Training Time/st} (seconds) 
& \textbf{Total Time} (mins)\\
\midrule
DeepONet & 1.08 & 0.10 & 86.3\\
FNO      & 0.57 & 18.54 & 30.9\\
\bottomrule
\end{tabular}

\vspace{1mm}
\begin{minipage}{0.92\linewidth}
\footnotesize
\textit{Note:} Training time is reported per step for DeepONet and per epoch for FNO.
\end{minipage}

\caption{Training parameters and time for DeepONet and FNO.}
\label{tab:runtime-3}
\end{table}

\subsection{Two-dimensional case (2V) - Boltzmann equation}

In this subsection, the numerical results for the Boltzmann equation in two dimensions are presented, where we adopt a point-to-point training strategy, i.e., the neural operator is trained between the initial distribution functions and the target distribution functions at a certain time. In addition, the modified FNO with conservative constraints is also tested in contrast to the standard FNO and the classical solver.
Although the present numerical experiments are restricted to spatially homogeneous problems, where the exact collision dynamics conserve the macroscopic moments, the learned neural operator does not automatically inherit these invariants. The conservative loss therefore acts as a physically motivated regularization mechanism that reduces moment drift in the surrogate model.

\subsubsection{Basic setup}

We establish a baseline configuration by fixing a standard set of model parameters as follows:
the number of Fourier modes $N_{F} = 16^2$, the number of hidden channels $N_{C} = 64$, the number of layers $N_{L} = 4$. In terms of the FNO with conservative constraints (C-FNO), the loss function consists of the standard mean squared error (MSE) and additional conservation constraints as in \eqref{C-Loss}.

The model is trained by the Adam optimizer with a learning rate of $10^{-3}$, a batch size of 8, and for 150 epochs. For the 2V Boltzmann equation, the experiments are conducted on a workstation equipped with an AMD Ryzen 9~9800X3D CPU, an NVIDIA GeForce RTX~5080 GPU, and 32\,GB of RAM.

As for the selection of the penalty coefficients $(\lambda_m,\lambda_p,\lambda_e)$ in \eqref{C-Loss}-\eqref{opt_cons}, in practice, we adopt an empirical but systematic calibration strategy. 
Firstly, we rescale the penalties so that the conservation terms $C_m, C_p, C_e$ are comparable in magnitude to the MSE term $\mathcal{L}_{\mathrm{MSE}}$ at the beginning of training, preventing either part from dominating the optimization. 
Then, to assess whether any single conservation law should be emphasized, we vary one coefficient while keeping the other two fixed, e.g., setting one of $(\lambda_m,\lambda_p,\lambda_e)$ to be one order of magnitude larger than the remaining two, and evaluate the resulting accuracy and conservation errors. 
In our experiments, no consistent advantage is observed by prioritizing a specific constraint; consequently, we use a unified weight $\lambda_m=\lambda_p=\lambda_e$ throughout the following numerical tests.
Finally, we tune the common value of $\lambda$ via repeated trials over a range of magnitudes and select the value that yields the best overall performance.

\subsubsection{Numerical experiment 5: BKW solution (Benchmark test)}

As one of the few analytical solutions available to the spatially homogeneous Boltzmann equation for Maxwell-type molecules, the exact BKW solution \eqref{BKW} serves as a benchmark for verifying the accuracy and stability of numerical algorithms.
In 2V, we can choose the benchmark values of the BKW solution \eqref{BKW} with $b(\cos\theta) = \frac{1}{2\pi}$, which then leads to \eqref{Kd}.

\begin{figure}[H]
  \centering
  \includegraphics[width=\linewidth]{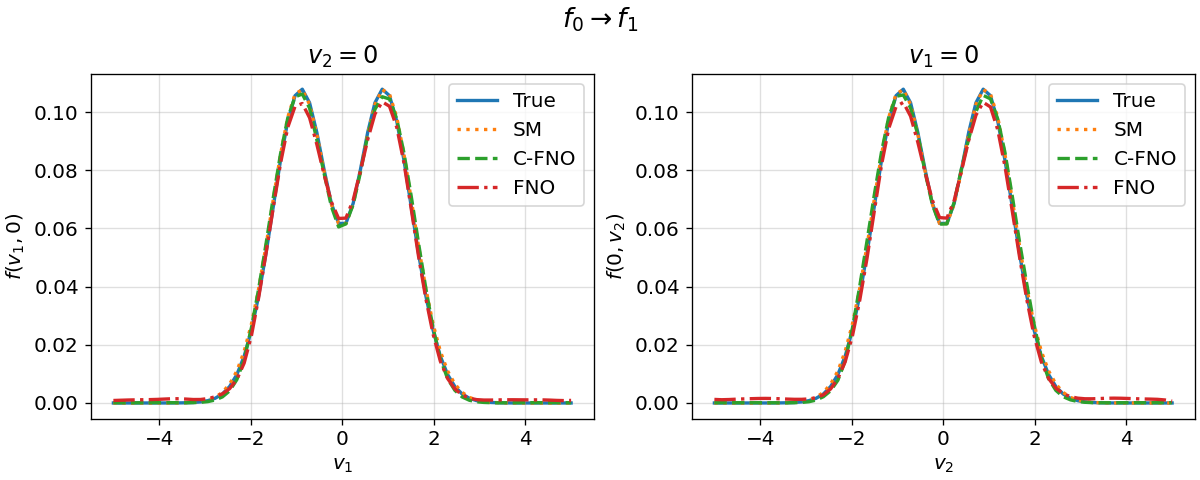}
  \caption{Profiles of the numerical solution $f^{\text{Num}}_{t=1}$ by different methods (SM, FNO, C-FNO) for $v_2{=}0$ (Left) and $v_1{=}0$ (Right) when $t=1$ (BKW solution).}
  \label{fig:bkw-t1}
\end{figure}
\begin{figure}[H]
  \centering
  \includegraphics[width=\linewidth]{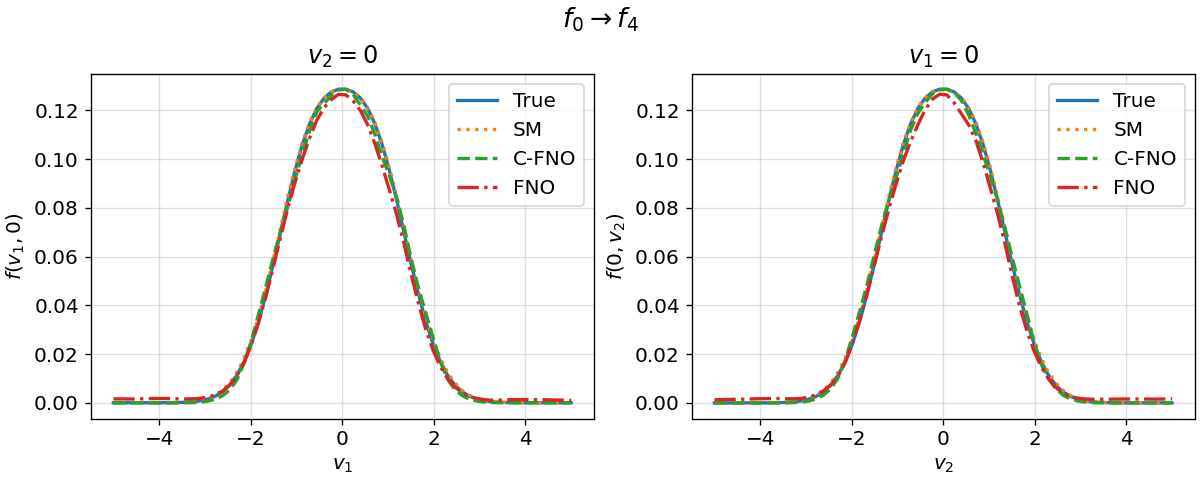}
  \caption{Profiles of the numerical solution $f^{\text{Num}}_{t=4}$ by different methods (SM, FNO, C-FNO) for $v_2{=}0$ (Left) and $v_1{=}0$ (Right) when $t=4$ (BKW solution).}
  \label{fig:bkw-t4}
\end{figure}
\begin{figure}[H]
  \centering
  \includegraphics[width=\linewidth]{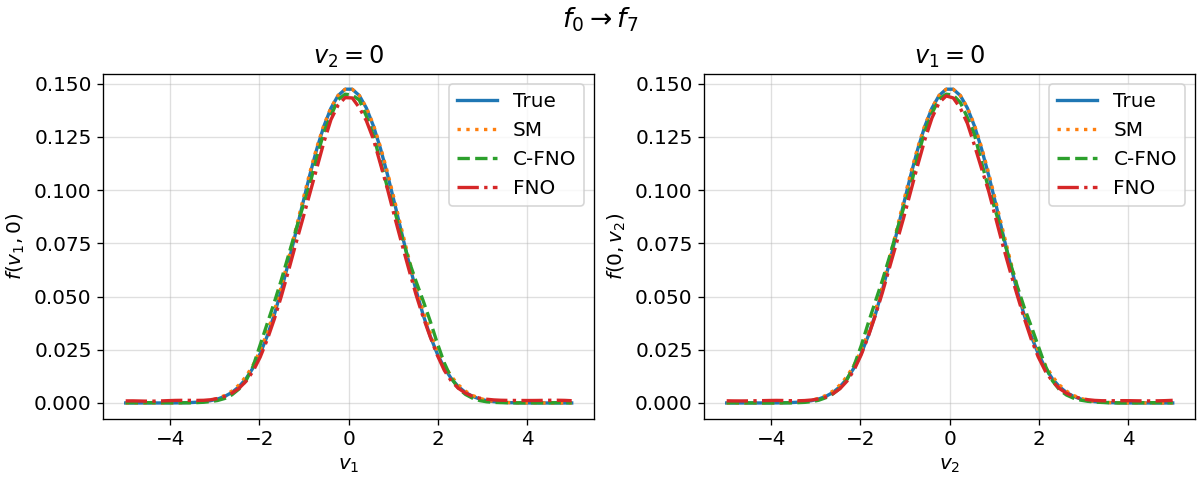}
  \caption{Profiles of the numerical solution $f^{\text{Num}}_{t=7}$ by different methods (SM, FNO, C-FNO) for $v_2{=}0$ (Left) and $v_1{=}0$ (Right) when $t=7$ (BKW solution).}
  \label{fig:bkw-t7}
\end{figure}

In what follows, we compare the profiles of the solutions obtained by the fast spectral method (SM), standard FNO framework (FNO) and our proposed FNO with conservative constraints (C-FNO) in short ($t=1$), intermediate ($t=4$) and long ($t=7$) time evolution, see Figure \ref{fig:bkw-t1}, \ref{fig:bkw-t4} and \ref{fig:bkw-t7}. Furthermore, the error of the solutions and their corresponding macroscopic quantities are presented in Table \ref{tab:bkw-t1}, \ref{tab:bkw-t2}, \ref{tab:bkw-t7} for different time scales, where the improvement of the conservative properties can be manifested when applying the constraints as in C-FNO.\\

\steptitle{(a) Case 5.1: Short-time evolution ($t=1$).}

\begin{table}[H]
\centering
\begin{tabular}{lccc}
\toprule
 Method & $L^1$ error & $L^2$ error & $L^{\infty}$ error \\
\midrule
 C-FNO  & 5.96e-02 & 1.09e-02 & 3.62e-03 \\
 FNO   & 9.04e-02 & 1.45e-02 & 5.14e-03 \\
 SM   & 4.11e-02 & 8.40e-03 & 2.57e-03 \\
\bottomrule
\end{tabular}
\begin{tabular}{ccc}
\toprule
  $| \rho^{\text{Num}} - \rho^{\text{True}} |$ & $\|u^{\text{Num}} - u^{\text{SM}}\|_2$ & $| E^{\text{Num}}-E^{\text{True}} |$ \\
\midrule
 2.11e-03 & 5.50e-03 & 2.68e-03 \\
 3.40e-02 & 2.06e-02 & 7.12e-01\\
 2.26e-13 & 6.18e-16 & 2.92e-06 \\
\bottomrule
\end{tabular}
\caption{$\|f^{\text{Num}}_{t=1}-f^{\text{True}}_{t=1}\|$ when $t=1$ (Left). Error in the macroscopic quantities (Right).}
\label{tab:bkw-t1}
\end{table}

\steptitle{(b) Case 5.2: Intermediate-time evolution ($t=4$).}

\begin{table}[H]
\centering
\begin{tabular}{lccc}
\toprule
 Method & $L^1$ error & $L^2$ error & $L^{\infty}$ error \\
\midrule
 C-FNO  & 8.03e-02 & 1.45e-02 & 4.88e-03 \\
 FNO   & 1.46e-01 & 2.49e-02 & 8.75e-03 \\
 SM   & 3.72e-02 & 7.20e-03 & 1.94e-03 \\
\bottomrule
\end{tabular}
\begin{tabular}{ccc}
\toprule
$|\rho^{\text{Num}} - \rho^{\text{True}}|$ & $\|u^{\text{Num}} - u^{\text{SM}}\|_2$ & $|E^{\text{Num}}-E^{\text{True}}|$ \\
\midrule
 2.56e-02 &8.51e-03  & 9.20e-02 \\
 2.75e-02 & 4.50e-02 & 8.87e-01\\
 2.26e-13 & 8.00e-16 & 7.53e-05\\
\bottomrule
\end{tabular}
\caption{$\|f^{\text{Num}}_{t=4}-f^{\text{True}}_{t=4}\|$ when $t=4$ (Left). Error in the macroscopic quantities (Right).}
\label{tab:bkw-t2}
\end{table}

\steptitle{(c) Case 5.3: Long-time evolution ($t=7$).}

\begin{table}[H]
\centering
\begin{tabular}{lccc}
\toprule
 Method & $L^1$ error & $L^2$ error & $L^{\infty}$ error \\
\midrule
 C-FNO  & 1.15e-01 & 2.18e-02 & 7.08e-03 \\
 FNO   & 1.52e-01 & 2.70e-02 & 9.02e-03 \\
 SM   & 3.85e-02 & 7.14e-03 & 1.83e-03 \\
\bottomrule
\end{tabular}
\begin{tabular}{ccc}
\toprule
$|\rho^{\text{Num}} - \rho^{\text{True}}|$ & $\|u^{\text{Num}} - u^{\text{SM}}\|_2$ & $|E^{\text{Num}}-E^{\text{True}}|$ \\
\midrule
 4.77e-04 & 3.23e-03 & 2.93e-02 \\
 3.67e-03 & 2.19e-02 & 5.64e-01\\
 2.26e-13 & 9.02e-16 & 3.09e-04 \\
\bottomrule
\end{tabular}
\caption{$\|f^{\text{Num}}_{t=7}-f^{\text{True}}_{t=7}\|$ when $t=7$ (Left). Error in the macroscopic quantities (Right).}
\label{tab:bkw-t7}
\end{table}

\subsubsection{Numerical experiment 6: Comparison of inference cost compared with the classical spectral solver}

We next compare the computational cost of the SM, FNO, and C-FNO at different resolutions. The presented times in Table~\ref{tab:time_comparison} correspond to the online prediction stage for one evaluation and do not include the offline training cost of the neural-operator models. This distinction is important: the neural-operator approach is most advantageous in many-query settings, where a trained model is repeatedly evaluated for different initial data or target times.

As shown in Table~\ref{tab:time_comparison}, the inference cost of FNO and C-FNO is lower than that of the spectral solver for the tested resolutions, and the gap becomes more visible as the grid size increases. Moreover, for point-to-point prediction from $f_{t=0}$ to $f_{t=1}$, the trained neural operator directly approximates the solution map, whereas the classical time-marching solver must advance through all intermediate time steps. Therefore, although training requires an offline cost, the learned operator can provide a computationally efficient surrogate once trained.

\begin{table}[htbp]
\centering
\begin{tabular}{c|ccc}
\hline
Resolutions $N$ & SM (seconds) & FNO (seconds) & C-FNO (seconds) \\
\hline
16  & $\approx$ 0.0162 & $\approx$ 0.0150 & $\approx$ 0.0151 \\
32  & $\approx$ 0.0595 & $\approx$ 0.0304 & $\approx$ 0.0301 \\
64  & $\approx$ 0.271 & $\approx$ 0.130 & $\approx$ 0.128\\
\hline
\end{tabular}
\caption{Comparison of inference computational time for one time step.}
\label{tab:time_comparison}
\end{table}

\subsubsection{Numerical experiment 7: Comparison with different training dataset sizes}

In this experiment, we further investigate how the size of the training dataset affects the accuracy of the proposed C-FNO method. The training samples are generated from the three classes of initial conditions described in Section~\ref{subsec:data-Boltzmann}. Thus, datasets of size \(300\), \(240\), and \(180\) correspond to \(100\), \(80\), and \(60\) samples, respectively, from each class of initial conditions. In all tests, the network architecture, optimizer, learning rate, and other training parameters are kept fixed.

\begin{figure}[H]
    \centering
    \includegraphics[width=\textwidth]{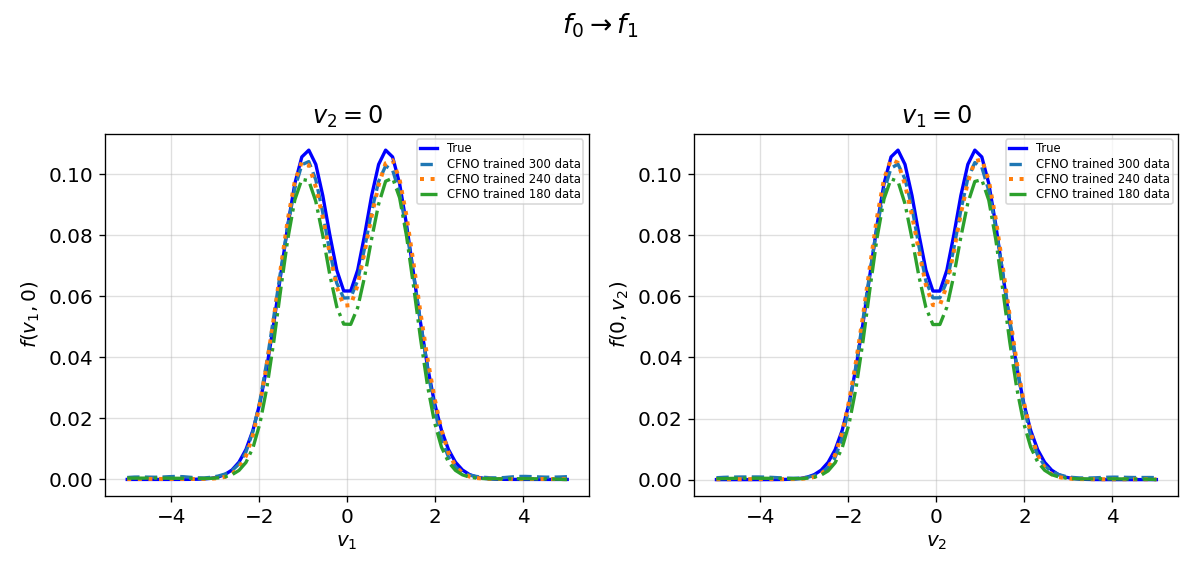}
    \caption{Profiles of the solution $f_{t=1}$ obtained by C-FNO with different numbers of training samples for $v_2{=}0$ (Left) and $v_1{=}0$ (Right).}
    \label{fig:data-size}
\end{figure}

Figure~\ref{fig:data-size} and Table~\ref{tab:dataset_error} show that reducing the number of training samples leads to a clear deterioration in accuracy. In particular, the \(L^2\)-error increases from \(1.84\times 10^{-2}\) for \(300\) samples to \(4.30\times 10^{-2}\) for \(240\) samples and \(1.27\times 10^{-1}\) for \(180\) samples. This trend indicates that the relatively large errors observed in some tests appear to be closely related to the limited amount of training data. Increasing the size and diversity of the training set is therefore expected to improve the accuracy of the learned operator, although larger-scale training is constrained by the computational resources available in the present study.

\begin{table}[H]
\centering
\begin{tabular}{c c c}
\hline
Number of training data & $L^2$ error & $L^{\infty}$ error \\
\hline
300 & 1.84e-02 & 4.78e-03 \\
240 & 4.30e-02 & 6.85e-03 \\
180 & 1.27e-01 & 1.37e-02\\
\hline
\end{tabular}
\caption{Error for different number of training data.}
\label{tab:dataset_error}
\end{table}

\subsubsection{Numerical experiment 8: General initial condition}

In addition to the benchmark test, we also generate a distinct class of initial conditions to test the generalization ability of the proposed framework.  
Specifically, we arbitrarily choose the initial condition as Gaussian type distribution with mean zero and a standard deviation randomly sampled from the range $[0.2, 1.0]$.  
This setup produces a variety of spread levels, ranging from sharply localized peaks to broadly dispersed distributions, thereby providing a diverse dataset for assessing long-time evolutions.

\begin{figure}[H]
  \centering
  \includegraphics[width=\linewidth]{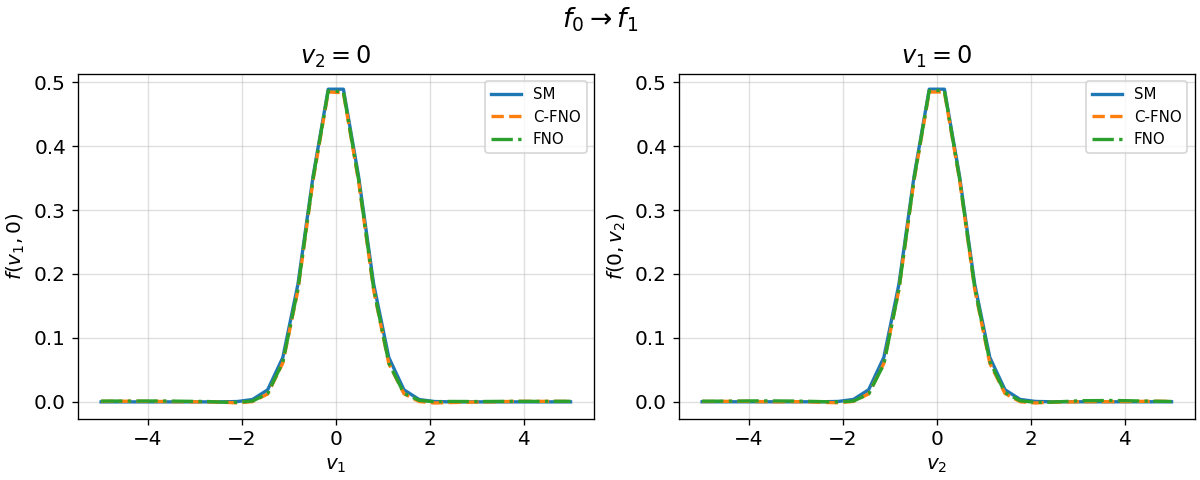}
  \caption{Profiles of the numerical solution $f^{\text{Num}}_{t=1}$ by different methods (SM, FNO, C-FNO) for $v_2{=}0$ (Left) and $v_1{=}0$ (Right) for $t{=}1$ (random initial condition).}
  \label{fig:gen-t1}
\end{figure}

\begin{figure}[H]
  \centering
  \includegraphics[width=\linewidth]{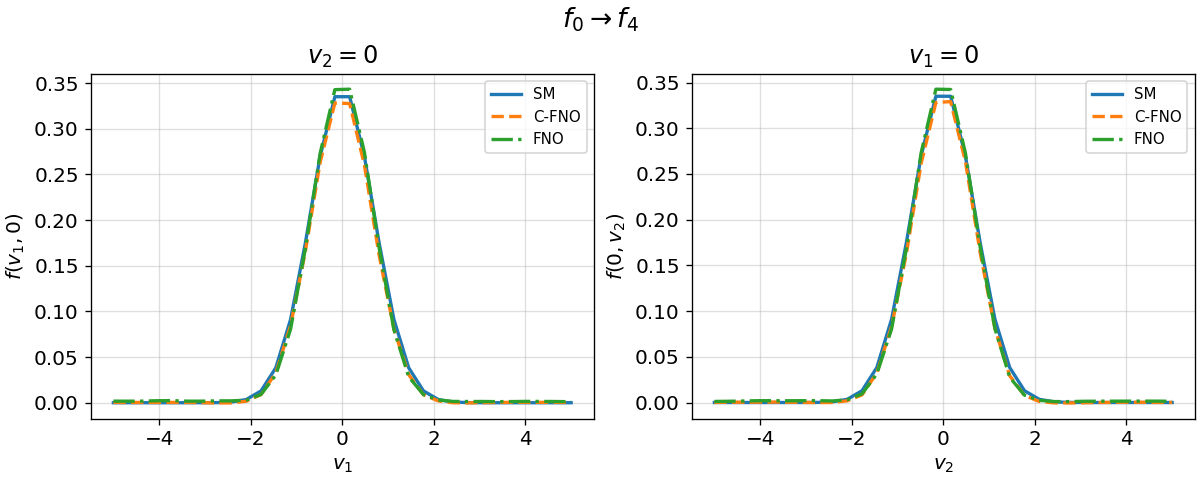}
  \caption{Profiles of the numerical solution $f^{\text{Num}}_{t=4}$ by different methods (SM, FNO, C-FNO) for $v_2=0$ (Left) and $v_1=0$ (Right) for $t=4$ (random initial condition).}
  \label{fig:gen-t4}
\end{figure}

\begin{figure}[H]
  \centering
  \includegraphics[width=\linewidth]{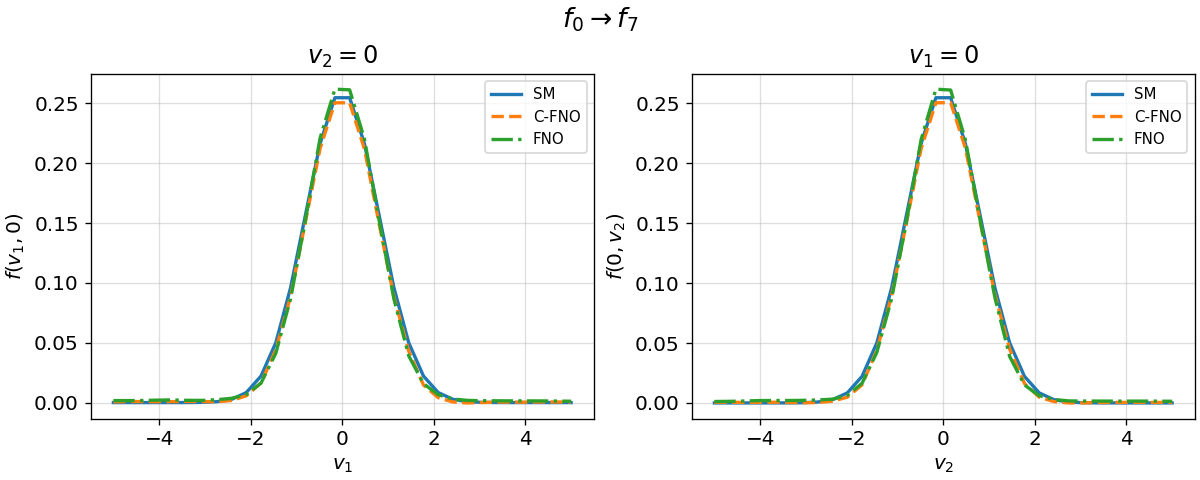}
  \caption{Profiles of the numerical solution $f^{\text{Num}}_{t=7}$ by different methods (SM, FNO, C-FNO) for $v_2=0$ (Left) and $v_1=0$ (Right) for $t=7$ (random initial condition).}
  \label{fig:gen-t7}
\end{figure}

In the following experiments, we compare the profiles of the solutions obtained by the standard FNO framework and our proposed C-FNO for general initial conditions, where the reference solutions are generated by the fast spectral method (SM) for comparison.
For the short ($t=1$), intermediate ($t=4$) and long ($t=7$) time evolution, we present the solution profiles in Figures \ref{fig:gen-t1}, \ref{fig:gen-t4} and \ref{fig:gen-t7}. To verify the strength of C-FNO in conservative property, the error of the solutions and their corresponding macroscopic quantities are presented in Table \ref{tab:gen-t1}, \ref{tab:gen-t4}, \ref{tab:gen-t7} for different time scales.\\

\steptitle{(a) Case 8.1: Short-time evolution ($t=1$).}

\begin{table}[H]
\centering
\begin{tabular}{lccc}
\toprule
 Method & $L^1$ error & $L^2$ error & $L^{\infty}$ error \\
\midrule
 C-FNO  & 7.73e-02 & 1.93e-02 & 8.55e-03 \\
 FNO   & 9.57e-02 & 2.49e-02 & 9.83e-03 \\
\bottomrule
\end{tabular}
\begin{tabular}{ccc}
\toprule
$|\rho^{\text{Num}} - \rho^{\text{SM}}|$ & $\|u^{\text{Num}} - u^{\text{SM}}\|_2$ & $|E^{\text{Num}}-E^{\text{SM}}|$ \\
\midrule
 2.02e-02 & 3.05e-03 & 1.25e-02 \\
3.29e-02 & 5.35e-03 & 1.86e-02\\
\bottomrule
\end{tabular}
\caption{$\|f^{\text{Num}}_{t=1}-f^{\text{SM}}_{t=1}\|$ when $t=1$ (Left). Error in the macroscopic quantities (Right).}
\label{tab:gen-t1}
\end{table}

\steptitle{(b) Case 8.2: Intermediate-time evolution ($t=4$).}

\begin{table}[H]
\centering
\begin{tabular}{lccc}
\toprule
 Method & $L^1$ error & $L^2$ error & $L^{\infty}$ error \\
\midrule
 C-FNO  & 9.71e-02 & 2.58e-02 & 1.38e-02 \\
 FNO   & 1.11e-01 & 3.30e-02 & 1.59e-02 \\
\bottomrule
\end{tabular}
\begin{tabular}{ccc}
\toprule
$|\rho^{\text{Num}} - \rho^{\text{SM}}|$ & $\|u^{\text{Num}} - u^{\text{SM}}\|_2$ & $|E^{\text{Num}}-E^{\text{SM}}|$ \\
\midrule
 1.16e-02 & 4.61e-03 & 7.85e-04 \\
2.41e-02& 6.62e-03 & 2.77e-02\\
\bottomrule
\end{tabular}
\caption{$\|f^{\text{Num}}_{t=4}-f^{\text{SM}}_{t=4}\|$ when $t=4$ (Left). Error in the macroscopic quantities (Right).}
\label{tab:gen-t4}
\end{table}

\steptitle{(c) Case 8.3: Long-time evolution ($t=7$).}

\begin{table}[H]
\centering
\begin{tabular}{lccc}
\toprule
 Method & $L^1$ error & $L^2$ error & $L^{\infty}$ error \\
\midrule
 C-FNO  & 1.08e-01 & 2.77e-02 & 1.28e-02 \\
 FNO   & 1.62e-01 & 4.60e-02 & 2.03e-02 \\
\bottomrule
\end{tabular}
\begin{tabular}{ccc}
\toprule
 $|\rho^{\text{Num}} - \rho^{\text{SM}}|$ & $\|u^{\text{Num}} - u^{\text{SM}}\|_2$ & $|E^{\text{Num}} - E^{\text{SM}}|$ \\
\midrule
 1.48e-02 & 2.44e-03 & 5.59e-03 \\
 2.54e-02 & 2.56e-03 & 2.81e-02\\
\bottomrule
\end{tabular}
\caption{$\|f^{\text{Num}}_{t=7}-f^{\text{SM}}_{t=7}\|$ when $t=7$ (Left). Error in the macroscopic quantities (Right).}
\label{tab:gen-t7}
\end{table}


\subsection{Three-dimensional case (3V)}
\label{sec:0d3v-bgk}

In this subsection, we present the numerical results for the spatially homogeneous BGK model and Boltzmann equation in a three-dimensional velocity space. We still adopt the sequence-to-sequence training strategy for the BGK model and point-to-point strategy for the Boltzmann equation. Compared with the 1V and 2V cases, the present 3V experiments are more demanding because the distribution function is represented on a higher dimensional velocity grid and the initial data contain strong multidirectional fluctuations.

\subsubsection{Numerical experiment 9: 3V BGK Model}
\label{subsec:3V-BGK}

We use the 3V BGK dataset constructed in Section~\ref{subsubsec:data-BGK}. The resulting $9900$ consecutive one-step pairs are randomly divided into $90\%$ training and $10\%$ validation sets.
For both FNO and C-FNO, the number of retained Fourier modes is $N_F=16$, the number of hidden channels is $N_C=32$, and the number of Fourier layers is $N_L=4$. The normalized distribution function and the three normalized velocity coordinates form the four input channels, while the distribution at the next time level is the single output channel. Both
models are optimized by the Adam optimizer for $50$ epochs with an initial learning rate of $10^{-3}$, a cosine annealing learning rate scheduler, and a batch size of two. The checkpoint with the smallest validation loss is retained for testing. The standard FNO is trained with the mean squared error loss. For C-FNO, the objective function additionally contains the mass, momentum, and energy conservation penalties defined in (3.13)--(3.14), with $\lambda_m=\lambda_p=\lambda_e=5\times10^{-3}$.
The training and evaluation are conducted on the same device as in Section \ref{subsec:1V-BGK}.

\begin{figure}[htbp]
    \centering
    \includegraphics[width=0.98\textwidth]{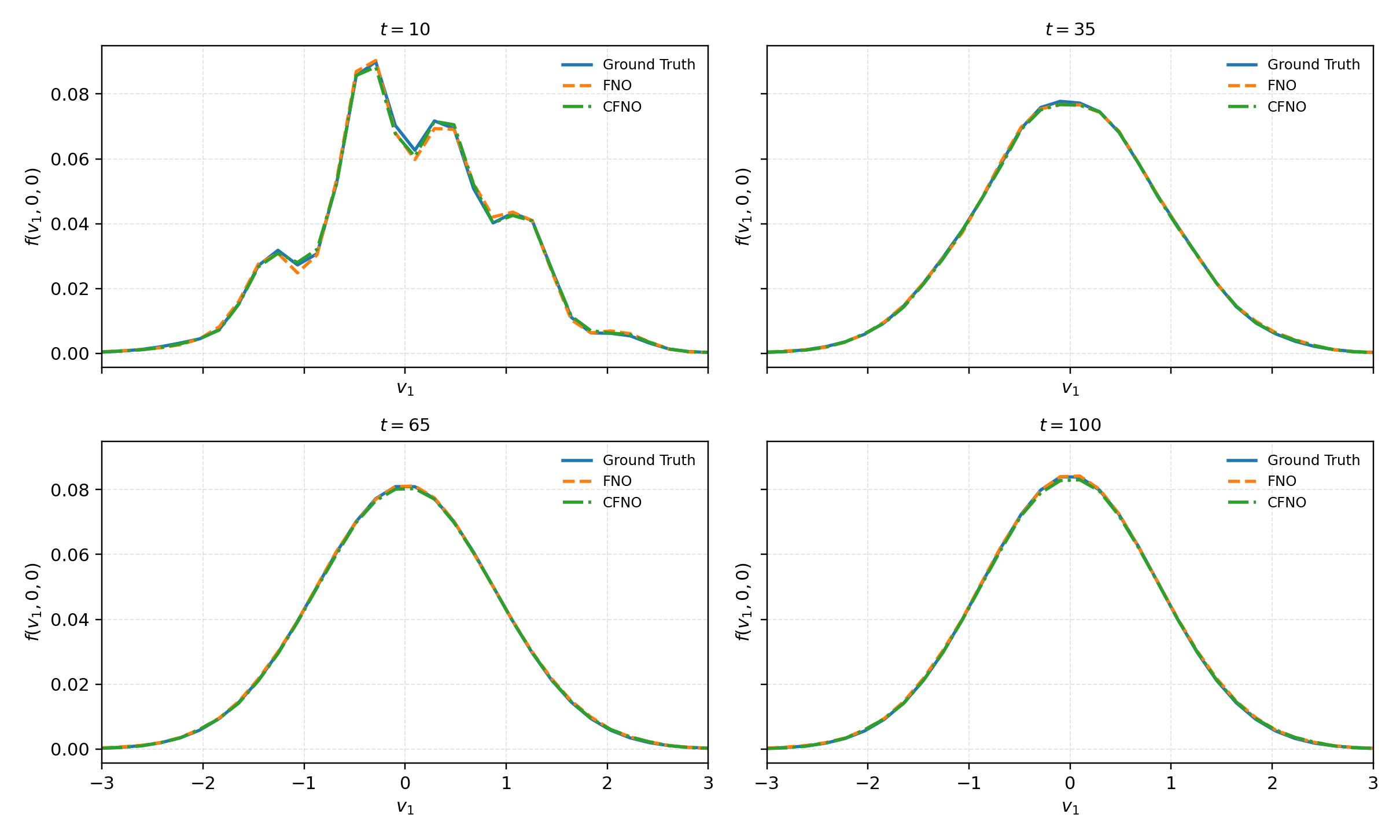}
    \caption{Profiles of the reference solution and numerical solution $f^{\text{Num}}_{t}(v_1, v_2 = 0, v_3 = 0)$ by different methods (FNO, C-FNO) at $t=10, 35, 65, 100$ for 3V BGK model.}
    \label{fig:0d3v-fno-cfno-profiles}
\end{figure}

\begin{figure}[htbp]
    \centering
    \includegraphics[width=0.78\textwidth]{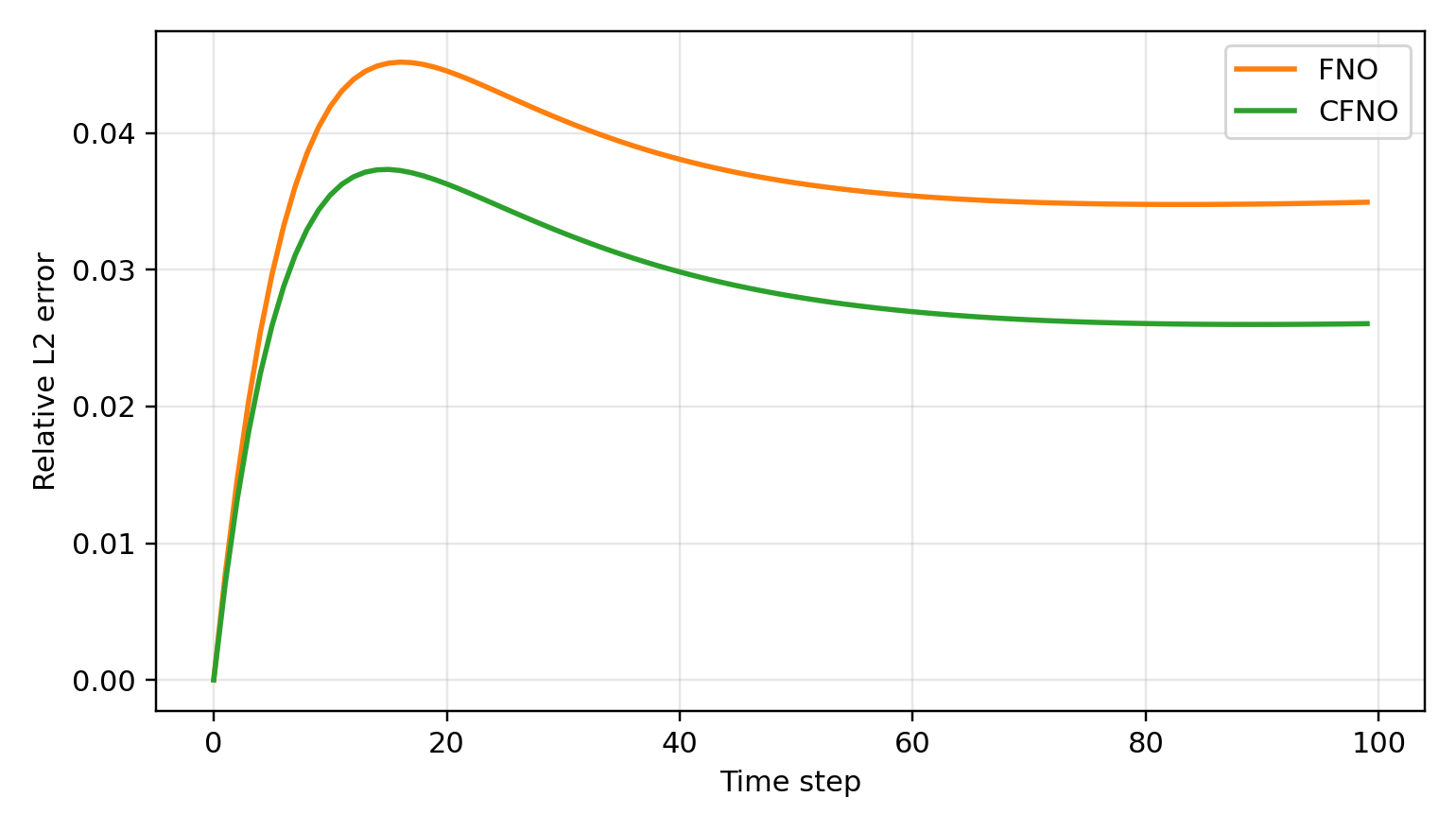}
    \caption{Evolution of relative $L^2$-error of FNO and C-FNO for 3V BGK model.}
    \label{fig:0d3v-fno-cfno-l2}
\end{figure}

Figure~\ref{fig:0d3v-fno-cfno-profiles} compares the reference solutions with the FNO and C-FNO numerical solution along the one-dimensional section $f(v_1,v_2=0,v_3=0)$ at different times. And we further quantify the accumulated prediction error by computing the $L^2$-error as Figure~\ref{fig:0d3v-fno-cfno-l2}. Both figures show that at the early stage, the distribution remains strongly oscillatory and far from equilibrium. Both models reproduce the principal peak and the smaller-scale structures accurately, while C-FNO is generally closer to the ground truth. As the BGK dynamics drive the distribution towards a Maxwellian, the three profiles become nearly indistinguishable. The agreement at the intermediate
and late stages indicates that both learned operators capture the dominant
relaxation dynamics under repeated application.

These results show that the standard FNO already provides an accurate and stable surrogate for the 3V BGK relaxation dynamics. Incorporating the conservation constraints further improves the one-step approximation and reduces error accumulation during autoregressive rollout. In particular, the C-FNO error remains consistently below the FNO error throughout the complete test horizon, demonstrating the benefit of physics-informed regularization for strongly fluctuating, high-dimensional distribution functions.

\subsubsection{Numerical experiment 10: 3V Boltzmann equation}
\label{subsec:3V-Boltzmann}

We next consider the spatially homogeneous Boltzmann equation for Maxwell
molecules in 3V velocity space.  As a benchmark, we choose the values of the BKW solution \eqref{BKW} with $d=3$.
Note that, unlike 2V case, the BKW solution \eqref{BKW} in 3V is nonnegative for
$t\geq 6\log(5/2)\approx5.498$. We therefore use the positive short time point-to-point ($t=5.5$ to $t=6.5$) as a benchmark test
\begin{equation}\label{eq:bkw-3V-benchmark-map}
  f_{t=5.5}(\cdot)
  \longmapsto f_{t=6.5}(\cdot).
\end{equation}

The velocity domain is discretized on a $ 32^3$ grid. The spectral method (SM) used the time step $\Delta t=0.1$, and the RK4 method in time discretization. Considering the training cost, both FNO and C-FNO are trained using the
same 50 one-unit trajectories with the number of Fourier modes $N_F=8^3$, hidden channels $N_{C} = 16$, layers $N_{L} = 4$. The exact pair in \eqref{eq:bkw-3V-benchmark-map} is reserved for evaluation and is not included in the training set.

\begin{figure}[htbp]
  \centering
  \includegraphics[width=0.94\textwidth]{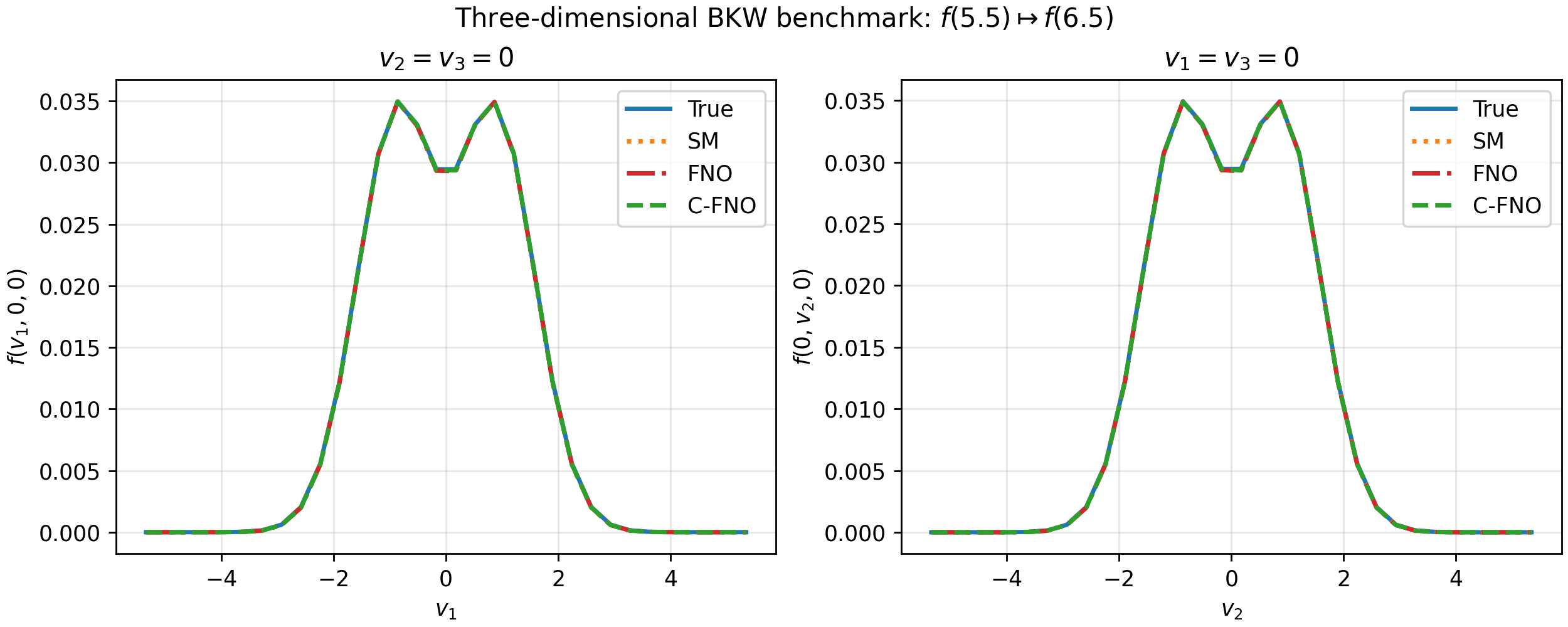}
  \caption{Profiles of the numerical solution $f^{\text{Num}}_{t=6.5}$ by different methods (SM, FNO, C-FNO) for $v_2 = v_3 = 0$ (Left) and $v_1 = v_3 = 0$ (Right) when $t=6.5$ (BKW solution).}
  \label{fig:bkw-3V-refined-profiles}
\end{figure}

In Figure~\ref{fig:bkw-3V-refined-profiles}, we compare the two cross profiles at $t=6.5$, where the profiles predicted by FNO and C-FNO nearly overlap the true BKW solution and the SM result. The corresponding errors are reported in Table~\ref{tab:bkw-3V-refined-errors}, from which one can find both FNO and C-FNO attain an error of the same order as the spectral reference for this benchmark test, and furthermore, C-FNO gives substantially smaller errors in the macroscopic quantities. 

\begin{table}[H]
\centering
\begin{tabular}{lccc}
\toprule
 Method & $L^1$ error & $L^2$ error & $L^{\infty}$ error \\
\midrule
 C-FNO & 8.5528e-03 & 6.7807e-04 & 2.3061e-04 \\
 FNO   & 8.2031e-03 & 6.9662e-04 & 2.3998e-04 \\
 SM    & 6.2891e-03 & 6.4930e-04 & 1.8569e-04 \\
\bottomrule
\end{tabular}
\begin{tabular}{ccc}
\toprule
 $|\rho^{\mathrm{Num}}-\rho^{\mathrm{True}}|$
 & $\|u^{\mathrm{Num}} - u^{\mathrm{True}}\|_2$
 & $|E^{\mathrm{Num}}-E^{\mathrm{True}}|$ \\
\midrule
 2.9350e-05 & 2.1609e-04 & 5.5386e-04 \\
 1.0602e-03 & 3.6495e-04 & 7.4118e-04 \\
 1.1939e-09 & 1.2517e-16 & 1.8704e-05 \\
\bottomrule
\end{tabular}
\caption{$\|f^{\text{Num}}_{t=6.5}-f^{\text{True}}_{t=6.5}\|$ when $t=6.5$ (Left). Error in the macroscopic quantities (Right).}
\label{tab:bkw-3V-refined-errors}
\end{table}

\subsubsection{Numerical experiment 11: Comparison of computing cost}
\label{subsec:numexp-3V-boltzmann-cost}

In this subsection, we compare the computing cost of the SM, FNO, and C-FNO for the 3V BKW benchmark. In Table~\ref{tab:timing-3V-boltzmann-online}, two different online SM times are informative: the first row is one RK4 time step of size $\Delta t=0.1$, while the second row is the complete time-marching solution from $t=5.5$ to $t=6.5$, consisting of ten RK4 steps. Since FNO and C-FNO learn the endpoint map directly, the fair online comparison with one neural inference is the complete SM evolution in ten time steps from $t=5.5$ to $t=6.5$ rather than a single SM time step.
All primary timings in Table~\ref{tab:timing-3V-boltzmann-online} were measured
on the same Apple M4 CPU for a batch of one on the $32^3$ velocity grid.  The SM was run in MATLAB double precision with precomputed spectral weights resident in memory.  Its reported times exclude weight precomputation and cache loading.  
Neural inference of FNO and C-FNO was run in PyTorch float32 with 20 warm-up evaluations and 200 synchronized repetitions.  
The neural operators reduce the online endpoint evaluation time by approximately three orders of magnitude.  
For reference, the one-time load of the $479.8$ MB spectral weight cache took
$2.09$s and is not included in the SM times above.
\begin{table}[htbp]
  \centering
  \small
  \begin{tabular}{llcc}
    \toprule
    Method & Online computing task & Time (seconds) & Speedup \\
    \midrule
    SM & One RK4 step, $\Delta t=0.1$
       & $ \approx 1.7490$ & -- \\
    SM & Complete map, $t=5.5\rightarrow6.5$
       & $\approx 17.5898$ & $1.0\times$ \\
    FNO & One endpoint inference, $t=5.5\rightarrow6.5$
       & $\approx 0.01429$ & $\approx 1231.0\times$ \\
    C-FNO & One endpoint inference, $t=5.5\rightarrow6.5$
       & $\approx 0.01440$ & $\approx 1221.5\times$ \\
    \bottomrule
  \end{tabular}
  \caption{Online computing time for 3V Boltzmann equation. Speedup is computed relative to the complete SM evolution in one complete step $t=5.5\rightarrow6.5$.}
  \label{tab:timing-3V-boltzmann-online}
\end{table}

Training is an offline expense and should not be combined with the cost of every subsequent prediction.  We therefore report it separately in Table~\ref{tab:timing-3V-boltzmann-training}.  Both models were trained for 3000 optimizer steps on the CPU using the same data split and deterministic
batch schedule.

\begin{table}[htbp]
  \centering
  \small
  \begin{tabular}{lccc}
    \toprule
    Method & Optimizer steps & Training time (seconds) & Training time (min) \\
    \midrule
    FNO   & 3000 & $\approx 416.8$ & $6.95$ \\
    C-FNO & 3000 & $\approx 417.3$ & $6.95$ \\
    \bottomrule
  \end{tabular}
  \caption{Offline training cost for FNO and C-FNO.}
  \label{tab:timing-3V-boltzmann-training}
\end{table}


\section{Conclusion}
\label{sec:conclusion}
In this work, we develop a Fourier Neural Operator (FNO)–based framework for learning the finite time evolution operator of the homogeneous Boltzmann equation and its simplified BGK model across multiple velocity space dimensions. The proposed methodology leverages the operator-learning paradigm to approximate the nonlinear mapping between distribution functions in various manners, where no fine-grained discretization is necessary. To enhance physical fidelity, conservation constraints were incorporated into the loss functional, ensuring better preservation of mass, momentum, and energy compared with the original FNO formulation. The framework presented in this paper can be naturally extended to higher-dimensional settings without requiring additional technical modifications.
Nevertheless, the use of soft constraints, i.e., enforcing conservation laws only through the loss function, remains insufficient to guarantee strict conservation, particularly when the training data are of limited quality. This limitation highlights the potential for our future work in the integration of hard physical constraints or structure-preserving correction mechanisms directly into the network architecture or training dynamics. 
Further extensions to more general collision kernels and spatially inhomogeneous kinetic equations also remain important directions for future study.
Such developments are expected to further strengthen the theoretical and computational foundation of physics-constrained operator learning, offering a promising new paradigm for kinetic theory and, more broadly, for nonlinear integro-differential equations arising in multiscale physical systems.

\section*{Acknowledgments} 
This work was initiated when KQ worked at the School of Mathematics, University of Minnesota - Twin Cities (UMN), and is part of BH's undergraduate research program at UMN under the mentorship of KQ. KQ thanks the kind support from the School of Mathematics of UMN during his working period, and also acknowledges support from AMS-Simons Travel Award grant. Part of this work is completed and based upon work supported by the National Science Foundation under Grant No.~DMS-2424139, while KQ was in residence at the Simons Laufer Mathematical Sciences Institute in Berkeley, California, during the Fall 2025 semester. The authors thank the two anonymous reviewers for their valuable comments and suggestions, which have significantly improved the quality of the manuscript.

\bibliographystyle{siam}
\bibliography{Qi_bib}

\end{document}